\documentclass{amsart}
\usepackage{amsmath,epsfig,amssymb,amsfonts,mathrsfs,url,color}

\newtheorem{thm}{Theorem}[section]
\newtheorem{lem}[thm]{Lemma}
\newtheorem{prop}[thm]{Proposition}
\newtheorem{cor}[thm]{Corollary}
\newtheorem{ex}[thm]{Example}
\newtheorem{defn}[thm]{Definition}
\numberwithin{equation}{section}

\def\Q{{\mathbb Q}}
\def\Z{{\mathbb Z}}
\def\C{{\mathbb C}}
\def\R{{\mathbb R}}

\def\sC{{\mathcal C}}
\def\sU{{\mathcal U}}
\def\sM{{\mathcal M}}
\def\sA{{\mathcal A}}

\def\u{{\bf u}}

\def\z{{\bf z}}
\def\y{{\bf y}}
\def\x{{\bf x}}
\def\dim{\mathrm{dim}}
\def\supp{\mathrm{supp}}
\def\arb {{\bf Arb}(\ell,n)}
\def\td{\operatorname{td}}
\def\vol{\operatorname{vol}}
\def\cone{\operatorname{cone}}

\begin{document}

\title[]{A Generating Function for all Semi-Magic Squares and 
the Volume of the Birkhoff Polytope}
\author{J.A. De Loera  \and   F. Liu  \and R. Yoshida }

\maketitle

\begin{abstract}
  We present a multivariate generating function for all $n \times
  n$ nonnegative integral matrices with all row and column sums equal
  to a positive integer $t$, the so called semi-magic squares. As a consequence we
  obtain formulas for all coefficients of the Ehrhart
  polynomial of the polytope $B_n$ of $n \times n$ doubly-stochastic
  matrices, also known as the Birkhoff polytope. In particular we derive
  formulas for the volumes of $B_n$ and any of its faces.
\end{abstract}

\section{Introduction}
\label{intro}

Let $B_n$ denote the convex polytope of $n \times n$ doubly-stochastic
matrices; that is, the set of real nonnegative matrices with all row
and column sums equal to one. The polytope $B_n$ is often called the
\emph{Birkhoff-von Neumann polytope},  the \emph{assignment polytope}, or
simply the \emph{Birkhoff polytope}. It is a well-known problem to compute
the volume of $B_n$ and there is a fair amount of work on the topic
(see \cite{beckpixton,chanrobbins,DG} and the references therein for
information on prior work); in this paper, we present  the first exact formula
for the volume of $B_n.$ The formula will follow from a multivariate
rational generating function for all possible $n \times n$ integer
nonnegative matrices with all row and column sums equal to a positive integer
$t$, the so called \emph{semi-magic squares} \cite{E,stanley1} (although
many authors refer to them as \emph{magic squares}).

Before stating our main formula, we give a few necessary definitions
and notation.  We call a directed spanning tree with all edges
pointing away from a root $\ell$ an \emph{$\ell$-arborescence}.  The
set of all $\ell$-arborescences on the nodes $[n]=\{1,2,\dots,n\}$
will be denoted by $\arb$. It is well known that the cardinality of
$\arb$ is $n^{n-2}$. For any $T \in \arb,$ we denote by $E(T)$ the set
of directed edges of $T.$ As usual let $S_n$ be the set of all
permutations on $[n].$ For any $\sigma \in S_n$, we associate $\sigma$ with
its corresponding permutation matrix, i.e., the $n \times n$ matrix
whose $(i,\sigma(i))$ entry is $1$ and zero otherwise. Throughout this
paper, we will use $\sigma$ to denote both a permutation and the
corresponding matrix and it should be clear which one it refers to
according to the context. The bracket operator $\langle \cdot, \cdot  \rangle$ denotes the dot product of two vectors.

It is well known that given a $d$-dimensional {\it integral polytope}
$P$, that is a polytope whose vertices have integer coordinates, for any
positive integer $t$, the number $e(P,t)$ of lattice points contained
in the $t$-th dilation, $t P = \{ tX \ | \ X \in P \}$, is a polynomial
of degree $\dim(P)$ in the variable $t$. Furthermore, the leading
coefficient of $e(P,t)$ is the \emph{normalized volume} of $P$ in units equal to
the volume of the fundamental domain of the affine lattice spanned by
$P$ (see Chapter 4 of \cite{stanley1} or the book
\cite{beckrobins}). This polynomial is called the {\it Ehrhart
  polynomial\/}\ of $P$.

One can find an expression for the Ehrhart polynomial
$e(B_n,t)$ of $B_n$ using the multivariate generating function
$$f(tB_n,\z)=\sum_{M \in tB_n \cap \Z^{n^2}} \z^{M}$$ of the
lattice points of $tB_n$, where $\z^{M} = \prod_{1 \le i, j \le n}
z_{i,j}^{m_{i,j}}$ if $M=(m_{i,j})$ is an $n$ by $n$ matrix in
$\R^{n^2}$. One can see that
by plugging $z_{i,j}=1$ for all $i$ and $j$ in $f(tB_n, \z),$ we get
$e(B_n,t).$ Our main result is

\begin{thm} \label{main}
Given any positive integer $t,$ 
the multivariate generating function for the lattice points of $tB_n$
is given by the expression
\begin{equation} \label{brionequ}
f(tB_n,\z)= \sum_{\sigma \in S_n} \sum_{T \in \arb} \z^{t\sigma}
\prod_{e \notin E(T)} \frac{1}{(1-\prod \z^{W^{T,e}\sigma})},
\end{equation}
where 
$\z^{t\sigma}=\prod_{k=1}^n z^t_{k,\sigma(k)}$.

Here $W^{T,e}$ denotes the $n \times n$ $(0,-1,1)$-matrix
associated to the unique oriented cycle in the graph $T+e$ (see
Definition \ref{dualrays} for details) and $W^{T,e} \sigma$ denotes the usual
matrix multiplication of $W^{T,e}$ and the permutation matrix
$\sigma$.  
\end{thm}

As we apply Lemma \ref{vollem} to Theorem \ref{main}, we
obtain the desired corollary:

\begin{cor} \label{volume} For any choice of fixed $\ell \in [n]$,
the coefficient of $t^k$ in the Ehrhart polynomial $e(B_n,t)$ of the
polytope $B_n$ of $n \times n$ doubly-stochastic matrices is given by
the formula
\begin{equation}
\frac{1}{k!}\sum_{\sigma \in S_n} \sum_{T \in \arb} \frac{{(\langle c, \sigma \rangle)^k}
\td_{d-k}(\{\langle c , W^{T,e}\sigma  \rangle, \ e \notin E(T) \} )}{\prod_{e \notin E(T)} \langle c , W^{T,e}\sigma \rangle}.
\end{equation}
In the formula $W^{T,e}$ is the $n \times n$ $(0,-1,1)$-matrix
associated to the unique oriented cycle in $T+e$ as defined in
Definition \ref{dualrays} and $W^{T,e} \sigma$ denotes the usual
matrix multiplication of $W^{T,e}$ and the permutation matrix
$\sigma$.  
The symbol $\td_{j}(S)$ is the $j$-th Todd polynomial
evaluated at the numbers in the set $S$ (see Definition \ref{todddef}
for details). Finally, $c \in \R^{n^2}$ is any vector such that
$\langle c ,W^{T,e}\sigma\rangle$ is non-zero for all pairs $(T,e)$ of
an $\ell$- arborescence $T$ and a directed edge $e \notin E(T)$ and all $\sigma
\in S_n$.

As a special case, the normalized volume of $B_n$ is given by 
\begin{equation}
\vol(B_n)= \frac{1}{((n-1)^2)!} \sum_{\sigma \in S_n} \sum_{T \in
\arb} \frac{\langle c, \sigma \rangle^{(n-1)^2}}{\prod_{e \notin
E(T)} \langle c, W^{T,e} \sigma\rangle}.
\end{equation}
\end{cor}

We stress that each rational function summand of Formula
\eqref{brionequ} is given only in terms of trees and cycles of a
directed complete graph. Our proof of Theorem \ref{main} is based in
the lattice point rational functions as developed in \cite{BarviPom}
with some help from the theory of Gr\"obner bases of toric ideals as
outlined in \cite{sturmfels}.

There is a large collection of prior work on this topic that we mention now
to put our result in perspective. In \cite{beckpixton} the authors
computed the exact value of the volume and the Ehrhart polynomials for
up to $n=10$, which is the current record for exact computation. The
computations in \cite{beckpixton} took several years of computer CPU
(running in a parallel machine setup) and our volume formula is so far
unable to beat their record without a much more sophisticated
implementation. On the other hand, in two recent papers, Canfield and
McKay \cite{canfieldmckay1,canfieldmckay2} provide simple asymptotic
formulas for the volume of $B_n$ as well as the number of lattice points
of $tB_n$. However, our closed formula for the volume of $B_n$ is nonetheless interesting for the following reasons. First, as it was demonstrated in
\cite{chanrobbinsyuen, BDV}, the faces of $B_n$ are also quite
interesting for combinatorics and applications. For example all
network polytopes appear as faces of a large enough $B_n$. From our
formula it is easy to work out volume formulas for \emph{any} concrete
face of $B_n$. We demonstrate this possibility in the case of the
well-known $CRY_n$ polytope \cite{chanrobbinsyuen} whose volume is
equal to the product of the first $n-1$ Catalan numbers (see \cite{zeilberger}).  Concretely,
we obtain for the first time the Ehrhart polynomials of facets of $B_n$ and
$CRY_n$ for  $n\leq 7$.  In principle, this could be applied to derive
formulas for the number of integral flows on networks.  Second, not only we can
derive formulas for the coefficients of the Ehrhart polynomial of
$B_n$, but we can also derive formulas for the integral of \emph{any}
polynomial function over $B_n$. We hope our generating function will
be useful for various problems over the set of all semi-magic squares, at least
for small values of $n$.

This paper is organized as follows: In Section \ref{ourmt} we begin
with background material that will be used in forthcoming sections,
including background  properties of $B_n$, a short discussion of Gr\"obner bases and
triangulations, Brion's theorem and generating functions for lattice points in polyhedra.
In that section, we sketch the steps we will follow to compute the generating function of lattice
points inside cones. In Section \ref{dualscones} we discuss the triangulations of
the dual cone at each vertex of $B_n$ which we encode via
Gr\"obner bases. From Brion's formula we derive in Section
\ref{formulalatticepts} a sum of rational functions encoding all the
lattice points of the dilation $tB_n$ and thus a proof of Theorem
\ref{brionequ}. In Section \ref{coeffsehrhart} we show how from
Theorem \ref{brionequ} we can derive all the coefficients of the
Ehrhart polynomial of $B_n$ after expressing the generating function
in terms of Todd polynomials.  Finally, in Section \ref{facesofBn}, we
explain how to obtain Ehrhart polynomials and formulas of integration 
for any face of $B_n$.

\section{Background}
\label{ourmt}

For basic definitions about
convex polytopes which are not stated in this paper, please see
\cite{ziegler}. Chapters 5 and 6 in \cite{YKK} have a very detailed
introduction to $B_n$ and transportation polytopes. For all the
details and proofs about lattice point counting and their multivariate
generating functions see \cite{barvibook,BarviPom,beckrobins}.
We begin with some useful facts about the polytope $B_n$. It
is well known that the vertices of $B_n$ are precisely the $n \times
n$ permutation matrices. Permutation matrices are in bijection with
matchings on the complete bipartite graph $K_{n,n}$. The polytope $B_n$ lies in the $n^2$-dimensional
real space $\R^{n^2} = \{ n \times n \text{ real matrices}\},$ and
we use $M(i,j)$ to denote the $(i,j)$-entry of a matrix $M$ in the space.
There is a graph theoretic description of the edges of $B_n$; they
correspond to the cycles in $K_{n,n}$. On the
other hand, for each pair $(i,j)$ with $1 \le i,j \le n$, the set of
doubly-stochastic matrices with $(i,j)$ entry equal to 0 is a {\it facet} (a
maximal proper face) of $B_n$ and all facets arise in this way. It is
also easy to see that the dimension of $B_n$ is $(n-1)^2$ (i.e., the
volume we wish to compute is the $(n-1)^2$-volume of $B_n$ regarded as
a subset of $n^2$-dimensional Euclidean space).  Note that an $n
\times n$ doubly-stochastic matrix is uniquely determined by its upper
left $(n-1)\times (n-1)$ submatrix.  The set of $(n-1)\times (n-1)$
matrices obtained this way is the set $A_n$ of all nonnegative
$(n-1)\times (n-1)$ matrices with row and column sums $\le 1$ such
that the sum of all the entries is at least $n-2$. $A_n$ is affinely
isomorphic to $B_n$ and we often compute in $A_n$ instead
of $B_n$ because $A_n$ is full-dimensional.

\subsection*{Cones and Generating functions for lattice points}

For any polytope $P \in \R^d,$ we would like to write a generating
function for the following sum encoding the lattice points of $P$
$$ \sum_{\alpha \in P \cap \Z^d} \z^{\alpha},$$
where $\z^{\alpha} = z_1^{\alpha_1} z_2^{\alpha_2} \cdots z_d^{\alpha_d}.$
We give now a step-by-step description of how the generating function 
is constructed.  

A \emph{cone} is the set of all linear nonnegative
combinations of a finite set of vectors. If a cone contains no other
linear subspace besides the origin then we say it is \emph{pointed}.
Given a cone $C \subset \R^d,$ the {\bf dual cone} to $C$ is a cone
$C^* = \{ y \in \R^d \ | \ \langle x, y \rangle \ge 0, \forall x \in
C\}.$ The following lemma states some properties of dual
cones. (See Theorem 9.1 in \cite{schrijver} for a proof).

\begin{lem}Let $C$ be a pointed cone in $\R^n,$  and let $D = C^*$ be 
its dual cone. Then the following properties hold
\begin{enumerate}\label{dualcone}
\item $C$ is the dual cone of $D$, namely $C = D^* = (C^*)^*.$
\item If $C$ is a full dimensional pointed cone, then so is
  D. Moreover, if $\{ F_i \}$ is the set of facets of $C,$ then $D$ is
  precisely the cone generated by the set of rays $\{ R_i \}$ satisfying, for
  any $i,$
\begin{equation}\label{dual}
R_i \mbox{ is perpendicular to } F_i, \mbox{ and for any ray $R$ of
$C$ not on $F_i:$ } \langle R, R_i \rangle > 0.
\end{equation}
\end{enumerate}
\end{lem}

Now if $P$ is a polytope and $v$ is a vertex of $P,$ the 
\emph{supporting polyhedron} of $P$ at $v$ is
$$S(P,v) = v + \{ u \in \R^d \ : \ v + \delta u \in P \ \mbox{for all
  sufficiently small } \delta > 0 \},$$ and the \emph{supporting cone}
of $P$ at $v$ is defined as $C(P,v) = S(P,v) - v.$

For a set $A \subset \R^d$, the indicator function $[ A ]: \R^d
\rightarrow \R$ of $A$ is defined as $$[ A ] (x) = \left \{
\begin{array}{ll} 1 \mbox{ if }x \in A, \\ 0 \mbox{ if }x \not \in
A.\\
\end{array}\right .$$ 
The algebra of polyhedra $\mathbb{P}(\R^d)$ is the vector space over
$\Q$ spanned by the indicator functions $[P]$ of all polyhedra $P
\subset \R^d$. The algebra of polytopes $\mathbb{P}_P(\R^d)$ is the
subspace spanned by the indicator functions of the polytopes in
$\R^d$. The algebra of cones $\mathbb{P}_C(\R^d)$ is the subspace
spanned by the indicator functions of the polyhedral cones in $\R^d$.
A linear transformation
\[
\Phi : \mathbb{P}(\R^d) \to V,
\]
where $V$ is a vector space over $\Q$ is called a {\em valuation}.  
Similarly, linear transformations defined on $\mathbb{P}_P(\R^d)$ and $\mathbb{P}_C(\R^d)$ are also called valuations \cite{BarviPom}.

One important tool for counting lattice points is the ability of
expressing the indicator function of a simplicial cone as an
integer linear combination of the indicator functions of unimodular simplicial
cones. Given a cone $K \subset \R^d$, we say
that the finite family of cones $K_i$, $i \in I = \{1,
2, \ldots , l\}$ is a \emph{decomposition} of $K$
if there are numbers $\epsilon_i \in \{ -1, 1 \}$ such that
$$[K] = \sum_{i \in I} \epsilon_i [K_i].$$

\begin{thm}[Theorem 3.1 and its proof in \cite{BarviPom}] \label{rational}
There is a map $\mathfrak F$ which, to each rational polyhedron $P
\subset \R^d$, associates a unique rational function $f(P,\z)$ in $d$ complex
variables $\z \in \C^d,$ $\z = (z_1, \dots, z_d),$ such that the
following properties are satisfied:
\begin{itemize}
\item[(i)] The map $\mathfrak F$ is a valuation. 
\item[(ii)] If $P$ is pointed, there exists a nonempty open subset
  $U_p \subset \C^d,$ such that $\sum_{\alpha \in P \cap \Z^d} \z^{\alpha}$ converges absolutely to
  $f(P, \z)$ for all $\z \in U_P.$ 
\item[(iii)] If $P$ is pointed,
then $f(P,\z)$ satisfies
$$f(P,\z) = \sum_{\alpha \in P \cap \Z^d} \z^{\alpha}$$ for any $\z
\in \C^d$ where the series converges absolutely.
\item[(iv)] If $P$ is not pointed, i.e., $P$ contains a line, then $f(P, \z) = 0.$ 
\end{itemize}
\end{thm}

Because the rational function $f(P,\z)$ encodes the lattice points of $P$,
we call $f(P,
\z)$ the {\it multivariate generating function of the lattice points}
(MGF) of $P$. The rational function has an expression as a sum of
simple terms, but to describe them we need the following facts.

\begin{thm}[Brion, 1988;  Lawrence, 1991](see \cite{BarviPom,beckhassesottile} 
for proofs) \label{brion} Let $P$ be a rational polyhedron and let $V(P)$
be the vertex set of $P$. Then,
$$f(P,\z) = \sum_{v \in V(P)} f(S(P,v),\z).$$
\end{thm}

This theorem reduces the problem of finding the MGF of a rational
polyhedron $P$ to that of finding the MGF of the supporting polyhedra at
each vertex of $P$. If the vertex of the supporting polyhedron is
integral we can simply assume the vertex is the origin and work instead with
supporting cones.  

\begin{cor}\label{intepoly} If $P$ an {\bf integral polyhedron}, i.e., 
all the vertices of $P$ are integral vertices, then
$$f(P,\z) = \sum_{v \in V(P)} \z^v f(C(P,v),\z).$$
\end{cor}

Although it is in general more complicated to give the MGF of an arbitrary
cone, if the cone is {\it unimodular}, its MGF has a simple form:

\begin{lem}[Lemma 4.1 in \cite{BarviPom}] \label{unicone}
If $K$ is a $d$-dimensional pointed cone in $\R^n$ generated by the rays $\{r_i\}_{1 \le i \le d},$ where
the $r_i$'s form a $\Z$-basis of the lattice $\Z^n \cap span(K)$ (span(K) is the $d$-subspace where $K$ lies), then we say $K$ is a {\bf unimodular cone} and we have that
$$f(K,\z) = \prod_{i=1}^d \frac{1}{1 - \z^{r_i}}.$$
\end{lem}

Barvinok gave an algorithm to decompose any pointed cone $C$ as a
signed sum of simple unimodular cones \cite{BarviPom} and thus
deriving an expression for $f(P,\z)$ as a sum of terms like those in
Lemma \ref{unicone}.  In principle, one needs to keep track of lower
dimensional cones in the decomposition for writing a
inclusion-exclusion formula of the MGF $f(C,\z).$ Fortunately, by
using the {\it Brion's polarization trick} (see Remark 4.3
in\cite{BarviPom}), one only needs to consider full-dimensional cones. 
This trick involves using dual cones of a decomposition of the dual cone to $C$
instead of directly decomposing $C$.  The main idea is to note  that  the duals of low dimensional cones
are not pointed and thus, from Part (iv) of Theorem \ref{rational}, their associated rational functions vanish.

Now we are ready to sketch the main steps of Barvinok's algorithm to compute $f(C,\z)$ (see \cite{BarviPom} for details): 
\\

\noindent {\bf Algorithm:} \\
Input: a rational cone full dimensional pointed cone $C.$ \\
Output: the MGF of $C: f(C, \z).$
\begin{enumerate}
\item Find the dual cone $D=C^*$ to $C$.
\item Apply the Barvinok decomposition to $D$ into 
a set of unimodular cones ${D_i}$ which have the same 
dimension as $D$ (ignoring all the lower dimensional cones).
\item Find the dual cone $C_i$ to each $D_i$.  The cone $C_i$ will be unimodular as well.
\item $f(C,\z) = \sum_{i} \epsilon_i f(C_i,\z),$ where $\epsilon_i$ is $+1$ or $-1$ determined by Barvinok decomposition.
\end{enumerate}

This algorithm is still not right for us; the algorithm is for full-dimensional cones, however, the cones we need to study
are not full-dimensional since the Birkhoff polytope is not full-dimensional. Also, Lemma \ref{dualcone}
provides us a way to compute the rays of $C^*$ if $C$
is full dimensional and pointed. Hence, it will be nice if we can make our cones full-dimensional.
What we will do is to properly project cones into a lower dimensional space so that 
they become full-dimensional. 

\begin{defn}\label{good}
Let $V \subset \R^n$ and $W \subset \R^m$ be vector spaces with full
rank lattices $L_V:= V \cap \Z^n$ and $L_W := W \cap \Z^m$,
respectively.  A linear map $\phi$ from $V$ to $W$ is a {\bf good
  projection} if $\phi$ gives a bijection between $L_V$ and $L_W.$
Note that because of the linearity of $\phi$,  the lattices
$L_V$ and $L_W$ have the same rank.
\end{defn}

\begin{lem}\label{goodform}
Suppose $V, W$ are as in Definition \ref{good} and $\phi$ is
a good projection from $V$ to $W.$ Because $\phi$ is a linear map, we
can consider $\phi$ is given by a certain $m \times n$ matrix $\phi =
(\phi_{i,j}).$ We define a map $\Phi: \C^m \to \C^n$ by mapping
$\y = (y_1, \dots, y_m) \in \C^m$ to $\z = (z_1, \dots, z_n) \in
\C^n,$ where
$$z_j = \prod_{i=1}^m y_i^{\phi_{i,j}}.$$
Then the following statements hold.
\begin{enumerate}
\item $\dim(V) = \dim(W).$
\item $\phi$ gives an isomorphism between $V$ and $W$ which preserves the lattice. Therefore, there exists an inverse (linear) map $\phi^{-1}$ from $W$ to $V$ that preserves the lattice as well. Thus, $\phi^{-1}$ is also a good projection from $W$ to $V.$
\item $C$ is a unimodular cone in $V$ if and only if $\phi(C)$ is a unimodular cone in $W.$
\item For any $\alpha \in \Z^n$ and $\y \in \C^m,$  if $\beta = \phi(\alpha)$ and $\z = \Phi(\y),$ then $\y^\beta = \z^\alpha.$
\item For any pointed rational polyhedron $P \in V,$ the series $\sum_{\beta \in \phi(P) \cap \Z^m} \y^\beta$ converges absolutely if and only if the series $\sum_{\alpha \in P \cap \Z^n} \Phi(\y)^\alpha$ converges absolutely. 
Furthermore, we have 
\begin{equation}\label{goodproj} f(\phi(P), \y) = f(P, \Phi(\y)).\end{equation}
\item Let $P_1, P_2, \dots, P_k$ be pointed rational polyhedra in $V,$ and $a_1, \dots, a_k \in \C,$ then 
$$f(P, \z) = \sum a_i f(P_i, \z) \Leftrightarrow f(\phi(P), \y) = \sum a_i f(\phi(P_i), \y).$$
\end{enumerate}
\end{lem}

\begin{proof}
The proofs of (1), (2) and (3) follow from the fact that good projections give lattices of same rank and thus isomorphic vector spaces. 
For the proof of (4), $\beta = \phi(\alpha)$ implies that $\beta_i = \sum_{j=1}^n \phi_{i,j} \alpha_j.$ Thus,
$$\z^\alpha = \prod_{j=1}^n z_j^{\alpha_j} = \prod_{j=1}^n \prod_{i=1}^m y_i^{\phi_{i,j} \alpha_j} 
= \prod_{i=1}^m \prod_{j=1}^n y_i^{\phi_{i,j} \alpha_j} = \prod_{i=1}^m  y_i^{\sum_{j=1}^n \phi_{i,j} \alpha_j} = \prod_{i=1}^m y_i^{\beta_i} = \y^\beta.$$

Because $\phi$ is a good projection, the lattice in $P$ and the lattice in $\phi(P)$ are in one-to-one correspondence under $\phi.$ Therefore, to prove (4), it is enough to show that if $\beta = \phi(\alpha)$ and $\z = \Phi(\y),$ then $\y^\beta = \z^\alpha.$ $\beta = \phi(\alpha)$ implies that $\beta_i = \sum_{j=1}^n \phi_{i,j} \alpha_j.$ Thus,
$$\z^\alpha = \prod_{j=1}^n z_j^{\alpha_j} = \prod_{j=1}^n \prod_{i=1}^m y_i^{\phi_{i,j} \alpha_j} 
= \prod_{i=1}^m \prod_{j=1}^n y_i^{\phi_{i,j} \alpha_j} = \prod_{i=1}^m  y_i^{\sum_{j=1}^n \phi_{i,j} \alpha_j} = \prod_{i=1}^m y_i^{\beta_i} = \y^\beta.$$

The first part of (5) follows immediately from (4). Let $Y$ be the set of $\y \in \C^m$ for which the series $\sum_{\beta \in \phi(P) \cap \Z^m} \y^\beta$ converges absolutely and  $Z$ be the set of $\z \in \C^n$ for which the series $\sum_{\alpha \in P \cap \Z^n} \z^\alpha$ converges absolutely. By the first part of (4), $\Phi(Y) \subset Z.$ By Theorem \ref{rational}, $f(P, \z) = \sum_{\alpha \in P \cap \Z^n} \z^\alpha$ for any $\z \in Z.$ In particular, $f(P, \z) = \sum_{\alpha \in P \cap \Z^n} \z^\alpha$ for any $\z \in \Phi(Y).$ Hence,  for any $\y \in Y.$
$$f(P, \Phi(\y)) = \sum_{\alpha \in P \cap \Z^n} \Phi(y)^\alpha = \sum_{\beta \in \phi(P) \cap \Z^m} \y^\beta.$$
We use Theorem \ref{rational} again to conclude that $f(P, \Phi(y))$ is the rational function $f(\phi(P), \y)$ associated to $\phi(P).$

Given (2), we only need to check one direction in (6). Suppose $f(P, \z) = \sum a_i f(P_i, \z).$ We can apply \eqref{goodproj} on both sides to obtain $ f(\phi(P), \y) = \sum a_i f(\phi(P_i), \y).$
\end{proof}

Using Lemma \ref{goodform}, we modify Barvinok's algorithm and sketch a method to construct $f(C,\z)$ for supporting cones $C$ at vertices
of $B_n$. We will try to follow this sequence of steps in  Section \ref{dualscones}:
\\

\noindent {\bf (CMGF) Method for constructing the multivariate generating function for lattice points of a cone:} \\
Input: a rational (not necessarily full dimensional)  pointed cone $C \subset \R^n.$ \\
Output: the MGF of $C:$ $f(C,\z)$.

\begin{enumerate}
\item[(0)] Let $V$ be the subspace spanned by $C$ in $\R^n.$ Find a subspace $W$ of $\R^m$ together with a good projection $\phi$ from $V$ to $W.$ Let $\overline{C} = \phi(C).$
\item Find a dual cone $\overline{D}$ to $\overline{C}$.
\item Decompose $\overline{D}$ into addition and subtraction of unimodular cones ${\overline{D}_i}$ which have the same dimension as $\overline{D}$, ignoring all the lower dimensional cones.
\item Find dual cone $\overline{C}_i$ of each $\overline{D_i}$. Note, that $\overline{C_i}$ is also unimodular. 
Let $C_i = \phi^{-1}(\overline{C}_i).$
\item $f(C,\z) = \sum_{i} \epsilon_i f(C_i,\z),$ where $\epsilon_i$ is $+1$ or $-1$ determined by the signed decomposition.
\end{enumerate}

In the next section, we will
apply the method (CMGF) step by step to the supporting cone at the
vertex $I$, the identity permutation. We will get the MGF of this
supporting cone and, by applying the action of the symmetric group
$S_n$, we can deduce the MGF of all other supporting cones of vertices
of $B_n$ and thus, by Theorem \ref{brion}, the MGF of $B_n.$
We will see later, in Section \ref{coeffsehrhart}, that the knowledge of
$f(P,\z)$ as a sum of rational functions yields a rational
function formula for the volume of $P$.

\subsection*{Triangulations and Gr\"obner bases of toric ideals}
For step $(2)$ in our step-by-step construction of the generating function, we will show (Lemma \ref{uni}) that in fact
any triangulation of the dual cone of the supporting cone of a vertex gives already a set of unimodular
cones (hence, the $\epsilon_i$'s in Step (4) are all $+1$). A triangulation of a cone $C$ is a special decomposition
of a cone as the union of simplicial cones with disjoint interiors whose union covers completely the cone $C$.
In this article we use polynomial ideals to codify the triangulations, namely  {\it
  toric ideals} and their {\it Gr\"obner bases}. See Chapter 8 in
\cite{sturmfels} for all details. Here are the essential notions:

Fix a set $\sA = \{a_1, a_2, \dots, a_n\} \subset \Z^d.$ For any $\u =
(u_1, u_2, \dots, u_n) \in \Z^n,$ we let
$$\u \sA : = u_1 a_1 + u_2 a_2 + \cdots + u_n a_n.$$ For any $\u \in
\Z^d,$ we denote by $\supp(\u) := \{ i \ | \ u_i \neq 0\}$ the {\it
support} of $\u.$ Every $\u \in \Z^d$ can be written uniquely as $\u =
\u^+ - \u^-,$ where $\u^+$ and $\u^-$ are nonnegative and have
disjoint support.

\begin{defn}
The toric ideal of $\sA$, $I_\sA \subset k[\x] := k[x_1, x_2, \dots, x_n]$ is the ideal generated
by the binomials
$$I_{\sA} := \langle \x^{\u^+} - \x^{\u^-} \ | \ \u \sA = 0 \rangle.$$
\end{defn}

     Given a real vector $\lambda=(\lambda_1,\dots,\lambda_n)$ in
     ${\mathbb R}^n$, we can define a \emph{monomial order}
     $>_\lambda$ that for any $a,b \in {\mathbb Z}_{\geq 0}^n$, their
     monomials satisfy $\x^a >_\lambda \x^b$ if $\langle a, \lambda
     \rangle > \langle b, \lambda \rangle$ and ties are broken via the
     lexicographic order. Using the ordering of monomials we can select the {\it initial monomial} of a polynomial $f$ with respect
     to $>_\lambda$, i.e., the highest term present. We will denote it by $in_{>_\lambda} (f)$. For an
     ideal $I$ contained in ${\mathbb C} \lbrack x_1,..,x_n \rbrack$
     its {\it initial ideal} is the ideal $in_{>_\lambda} (I)$
     generated by the initial monomials of all polynomials in $I$. A
     finite subset of polynomials $G=\lbrace g_1,...,g_n \rbrace$ of
     an ideal $I$ is a {\it Gr\"obner basis} of $I$ with respect to
     $>_\lambda$ if $in_{>_\lambda} (I)$ is generated by $\lbrace
     in_{>_\lambda} (g_1),...,in_{>_\lambda} (g_n) \rbrace$. In other
     words, $G$ is a Gr\"obner basis for $I$ if the initial monomial
     of any polynomial in $I$ is divisible by one of the monomials
     $in_{>_\lambda} (g_i)$. It can be proved from the definition that a
     Gr\"obner basis is a generating set for the ideal $I$.  As we
     will state later, each Gr\"obner basis of the toric ideal
     $I_{\sA}$ yields a \emph{regular triangulation} of the convex
     hull of $\sA$.  The fact that triangulations constructed using
     Gr\"obner bases are regular will not be used in our construction.

A {\em subdivision} of $\sA$ is a collection $T$ of subsets of $\sA$,
called {\em cells}, whose convex hulls form a polyhedral complex with
support $Q=conv(\sA)$. If each cell in $T$ is a simplex, then $T$ is
called a {\em triangulation} of $\sA$.  Every vector
$\lambda=(\lambda_1,\dots,\lambda_n)$ in ${\mathbb R}^n$ induces a
subdivision of $\sA=\{a_1,\dots,a_n \}$ as follows. Consider the
polytope $Q_\lambda=conv(\{(a_1,\lambda_1),\dots,(a_n,\lambda_n)\})$
which lies in ${\mathbb R}^{d+1}$. Generally, $Q_\lambda$ is a
polytope of dimension $dim(conv(\sA))+1$.  The {\em lower envelope} of
$Q_\lambda$ is the collection of faces of the form $ \{ x \in
Q_\lambda \vert \langle c, x \rangle =c_0\}$ with $Q_\lambda$ contained in the halfspace
$\langle c, x \rangle \leq c_0$ and the last coordinate $c_{d+1}$ is negative.  The
lower envelope of $Q_\lambda$ is a polyhedral complex of dimension
$dim(conv (\sA))$. We define $T_\lambda$ as the subdivision of
$\sA$ whose cells are the projections of the cells of the lower
envelope of $Q _\lambda$. In other words, $\{ a_{i_1},
a_{i_2},\dots,a_{i_k}\}$ is a cell of $T_\lambda$ if
$\{(a_{i_1},\lambda_{i_1}),(a_{i_2},\lambda_{i_2}),\dots,(a_{i_k},\lambda_{i_k})\}$
are the vertices of a face in the lower envelope of $Q_\lambda$.  The
subdivision $T_\lambda$ is called a \emph{regular subdivision} of
$\sA$. Remark that just as a triangulation can be uniquely specified
by its maximal dimensional simplices, it can also be uniquely
expressed by its minimal non-faces (minimal under containment).  Now
we are ready to state the algebra-triangulation correspondence:

\begin{thm}[See proof in Chapter 8 of \cite{sturmfels}] \label{bernds}
Let $\sA$ be an $n \times d$ matrix with integer entries, whose rows
vectors $\{a_1,\dots,a_n \}$ span an affine space of dimension
$d-1$. Let $I_\sA$ be the toric ideal defined by $\sA$. Then, the 
minimal non-faces of the regular triangulation of $\sA$ associated to the
vector $\lambda$ can be read from the generators of the radical of the 
initial ideal of the Gr\"obner basis of $I_\sA$ with respect to the
term order $>_{\lambda}$. More precisely, for $\lambda$ generic, the radical of the initial ideal
of $I_\sA$ equals
$$ \langle x_{i_1}x_{i_2}\cdots x_{i_s}: \{i_1,i_2,\dots,i_s\} \ \mbox{is a minimal non-face of $T_\lambda$} \rangle = \bigcap_{\sigma \in T_\lambda} 
\langle x_i : i \not\in \sigma \rangle. $$
                    
\end{thm}

The crucial fact we will use is that the maximal simplices of the regular triangulation $T_\lambda$ are
transversals to the supports of the monomials from the initial ideal of the Gr\"obner basis.
In the next section, we will apply Theorem \ref{bernds}
to create a triangulation of the dual cones. 

To the readers who are unfamiliar with commutative algebra language,  using
a Gr\"obner basis to describe a triangulation may not feel totally necessary or clear. Thus,
we explain here the advantages of doing it this way. First, traditionally checking that a set
of simplices is a triangulation of $\sA$ is not trivial since one has to verify they have disjoint
interiors (which requires a full description of all linear dependences of the
rays) and that the union of the simplicial cones fully covers the convex hull of $\sA$.
But, having a Gr\"obner basis avoids checking these two tedious geometric facts.
Second, the initial monomials of the Gr\"obner bases are precisely
the minimal non-faces of the triangulation $T_\lambda$, which are complementary to the maximal
simplicial cones of the triangulation. From the point of view of efficiency, the encoding
of a simplicial complex via its non-faces is sometimes much more economic than via its maximal facets. 
For more on the theory of triangulations see \cite{DRS}.

\section{The MGF of the supporting cone of $B_n$ at the vertex $I$} \label{dualscones}


Due to the transitive action of the symmetric group on $B_n$ it is enough
to explain a method to compute the MGF of the supporting cone at the vertex associated
to the identity permutation (we denote this by $I$)
and then simply permute the results. Nevertheless
it is important to stress that, although useful
and economical, there is no reason to use the
same triangulation at each vertex. Similarly, the
triangulations we use are all regular, but for our
purposes there is no need for this property either.

There are $n^2$ facets of $B_n:$ for any fixed $(i,j): 1 \le i, j \le n,$ the
collection of permutation matrices $P$ satisfying $P(i,j) = 0$ defines
a facet $F_{i,j}$ of $B_n.$ Hence, every permutation matrix is on
exactly $n(n-1)$ facets and the vertex $I$ is on the facets $F_{i,j},
i \neq j.$
Let $C_n$ be the supporting cone at the identity matrix $I,$
then the set of facets of $C_n$ is $\{F_{i,j} - I \}_{1 \le i, j \le n, i \neq j}.$ (Note that we need to subtract the vertex $I$ from $F_{i,j}$ because the supporting cone is obtained by shifting the supporting polyhedron at the vertex $I$ to the origin.) We are going to apply our method CMGF to find the MGF of $C_n.$

\subsection{Step 0: A good projection}
$C_n$, as well $B_n$, lie in the $n^2$-dimensional space $\R^{n^2} = \{ n \times n$ real matrices$\}.$ But they lie in different affine subspaces (the vertex of $C_n$ 
is the origin). Let $V_n$ be the subspace of $\R^{n^2}$ spanned by $C_n.$ It is easy to see that 
\begin{equation}\label{sub}
V_n = \{ M \in \R^{n^2} \ | \ \sum_{k=1}^n M(i,k) = \sum_{k=1}^n M(k,j) = 0, \forall i, j\}.
\end{equation} 
Let $W_n$ be the vector space $\R^{(n-1)^2} = \{ (n-1) \times (n-1)$ real matrices$\}.$ We define a linear map $\phi$ from $\R^{n^2}$ to $W_n$ by ignoring the entries in the last column and the last row of a matrix in $\R^{n^2},$ that is, for any $M,$ we define $\phi(M)$ to be the matrix $(M(i,j))_{1 \le i, j \le n-1}.$ 
One can check that when restrict $\phi$ from $V_n$ to $W_n$,  $\phi$ is a good projection from $V_n$ to $W_n$.  Let $$\overline{C}_n := \phi(C_n).$$ 

Also, let $\overline{F}_{i,j} = \phi(F_{i,j})$ and $\overline{P} = \phi(P),$ for any permutation matrix $P$ on $[n].$ (These are actually the facets and vertices for $A_n$ which is the full-dimensional version of $B_n$ we explained at the beginning of Section \ref{ourmt}.) By the linearity of $\phi$
the facets of $\overline{C}_n$ are $\{\overline{F}_{i,j} - \overline{I} \}_{1 \le i, j \le n, i \neq j},$ and $\overline{F}_{i,j}$ is defined by the collection of $\overline{P}$'s where $P$'s are permutation matrices (on $[n]$) satisfying $P(i,j) = 0.$

\subsection{Step 1: The dual cone $\overline{D}_n$ to $\overline{C}_n$}
The cone $\overline{C}_n$ is full dimensional in $W = \R^{(n-1)^2}.$ Hence, we can use Lemma \ref{dualcone} to find its dual cone. 
We will first define a cone, and then show it is the dual cone to $\overline{C}_n.$





\begin{defn}
$\overline{D}_n$ is the cone spanned by rays
$\{\overline{M}_{i,j}\}_{1\le i,j \le n, i \neq j},$ where
$\overline{M}_{i,j}$ is the $(n-1)$ by $(n-1)$ matrix such that

\begin{itemize}
\item[(i)] the $(i,j)$-entry is $1$ and all other entries equal zero, if $i \neq n$ and $j \neq n;$
\item[(ii)] the entries on the $i$th row are all $-1$ and all other entries equal zero, if $i \neq n$ and $j = n;$
\item[(iii)] the entries on the $j$th column are all $-1$ and all other entries equal zero, if $i = n$ and $j \neq n.$
\end{itemize}

\end{defn}

\begin{ex}[Example of $\overline{M}_{i,j}$ when $n = 3$] \label{exM1} Here we present each $2$ by $2$ matrix $\overline{M}_{i,j}$ as a row vector, which is just the first and second row of the matrix in order.

$\begin{array}{cccccc} \overline{M}_{1,3}: & -1 & -1 & & 0 & 0 \\
\overline{M}_{2,3}: & 0 & 0 & & -1 & -1 \\ \overline{M}_{3,1}: & -1 &
0 & & -1 & 0 \\ \overline{M}_{3,2}: & 0 & -1 & & 0 & -1 \\
\overline{M}_{1,2}: & 0 & 1 & & 0 & 0 \\ \overline{M}_{2,1}: & 0 & 0 &
& 1 & 0
\end{array}$
\end{ex}

\begin{lem}\label{dualkn}
$\overline{D}_n$ is the dual cone to $\overline{C}_n$ inside the vector space $W_n$. 
\end{lem}

\begin{proof}
For any $i, j \in [n]$ and $i \neq j,$ we need to check that condition
(\ref{dual}) is satisfied. Note that a ray of $\overline{C}_n$ is given by the vector $\overline{P}-\overline{I}$, for $P$ a permutation
matrix adjacent to the identity permutation. Thus it is enough to show that for any
permutation matrix $P$ on $[n],$ we have $\langle \overline{M}_{i,j}, \overline{P} \rangle
\ge \langle \overline{M}_{i,j}, \overline{I} \rangle$ and the equality holds if and only if
$\overline{P}$ is on $\overline{F}_{i,j}$, or equivalently, $P$ is on the facet $F_{i,j}.$ We have the following three situations for verification:

\begin{itemize}
\item[(i)] If $i \neq n$ and $j \neq n,$ $\langle \overline{M}_{i,j}, \overline{P} \rangle$ is $0$ if $P$ is on $F_{i,j}$ and is $1$ if $P$ is not on $F_{i,j}.$
\item[(ii)] If $i \neq n$ and $j = n,$ $\langle \overline{M}_{i,j}, \overline{P} \rangle$ is $-1$ if $P$ is on $F_{i,j}$ and is $0$ if $P$ is not on $F_{i,j}.$
\item[(iii)] If $i = n$ and $j \neq n,$ it is the same as (ii).
\end{itemize}
Therefore, $\overline{D}_n$ is the dual cone to $\overline{C}_n.$ 
\end{proof}

\subsection{Step 2: The triangulations of $\overline{D}_n$}
As we mentioned in the last section, we will use the idea of toric ideal to find a triangulation of the dual cone $\overline{D}_n$ to decompose $\overline{D}_n$ into unimodular cones.

\begin{lem}\label{uni} Let
$\sM$ be the configuration of vectors $\{\overline{M}_{i,j}\}_{1\le
i,j \le n, i \neq j}$ and $[\sM]$ denote the matrix associated to
$\sM$, i.e, the rows of $[\sM]$ are the vectors in $\sM$ written as row vectors. The matrix $[\sM]$ is totally unimodular, i.e., for any
$(n-1)^2$ linearly independent $\overline{M}_{i,j}$'s, they span a
unimodular cone. It follows that all triangulations of the cone $\overline{D}_n$
have the same number of maximal dimensional simplices.
\end{lem}

\begin{proof}
Up to a rearrangement of rows the matrix $[\sM]$ will look as follows: 
The first few rows are the negatives of the vertex-edge incidence matrix of the
complete bipartite $K_{n-1,n-1}$, then under those rows we have $n-1$
cyclically arranged copies of an $(n-2) \times (n-2)$ identity matrix. It
is well known that the vertex-edge incidence matrix of the complete
bipartite $K_{n-1,n-1}$ is totally unimodular. Moreover it is also known,
see e.g., Theorem 19.3 in \cite{schrijver}, that a matrix $A$ is
totally unimodular if each collection of columns of $A$ can be split
into two parts so that the sum of the columns in one part minus the
sum of the columns in the other part is a vector with entries only
$0,+1$, and $-1$. This characterization of totally unimodular matrices
is easy to verify in our matrix $[\sM]$ because whatever partition that
works for the columns sets of the vertex-edge incidence matrix of the
complete bipartite $K_{n-1,n-1}$ works also for the corresponding
columns of $\sM$, because the diagonal structure of the rows below it.

The fact that all triangulations have the same number of maximal simplices
follows from the unimodularity as proved in Corollary 8.9 of \cite{sturmfels}.
\end{proof}

Therefore, any triangulation of $\overline{D}_n$ gives a decomposition of 
$\overline{D}_n$ into a set of unimodular cones. Since $\sM$ defines
the vertex figure of $\overline{D}_n,$ it is sufficient to triangulate
the convex hull of $\sM.$ Hence, we consider the toric ideal
$$I_{\sM} := \langle \x^{\u^+} - \x^{\u^-} \ | \ \u \sM = 0 \rangle$$
of $\sM$ inside the polynomial ring $k[\x] := k[x_{i,j} : 1\le i,j\le
n, i \neq j].$ Note here $\u \in \Z^{n(n-1)}$ is a $n(n-1)$
dimensional vector indexed by $\{(i,j): i,j \in [n], i \neq j\}.$

Recall that a {\it circuit} of $I_{\sA}$ is an irreducible binomial
$\x^{\u^+} - \x^{\u^-}$ in $I_{\sA}$ which has minimal
support. Another result follows immediately from Lemma \ref{dualkn},
Lemma \ref{uni} and \cite[Proposition 4.11, Proposition 8.11]{sturmfels}:

\begin{lem}
The set $\sC_{\sM}$ of circuits of the homogeneous toric ideal $I_{\sM}$  
is in fact a universal Gr\"obner basis $\sU_{\sM}$ for $I_{\sM}$.
\end{lem}

For any partition of $[n] = S \cup T,$ we denote by $\u_{S,T} \in \Z^{n(n-1)}$ the $n(n-1)$ dimensional vector, where
$$\u_{S,T}(i,j) = \begin{cases}1, &  \mbox{if $i \in S, j \in T$,}  \\
        -1,  &  \mbox{if $i \in T, j \in S$},\\
        0, &  \mbox{otherwise}.\end{cases}$$
One can easily check that $\u_{S,T}$ has the following two properties:
\begin{equation}\label{anti}
\u_{S,T}(i,j) + \u_{S,T}(j,i)=0, \mbox{ for any $i \neq j$.}
\end{equation}
\begin{equation}\label{jacob}
\u_{S,T}(i,j) + \u_{S,T}(j,k)+ \u_{S,T}(k,i) = 0, \mbox{ for any distinct $i, j$ and $k.$}
\end{equation}

We define
$$P_{S,T} := \x^{\u_{S,T}^+} - \x^{\u_{S,T}^-} =\prod_{s \in S, t \in T} x_{s,t} - \prod_{s \in S, t \in T} x_{t,s},$$
where $\u_{S,T}^+(i,j) = \begin{cases}1, &  \mbox{if $i \in S, j \in T$},  \\
        0, &  \mbox{otherwise},\end{cases}$
and $\u_{S,T}^-(i,j) = \begin{cases}1, &  \mbox{if $i \in T, j \in S$} , \\
        0, &  \mbox{otherwise}.\end{cases}$

\begin{prop}\label{circ} The set of circuits of $I_{\sM}$ consists of all the binomials $P_{S,T}$'s:
$$\sC_{\sM} = \{ P_{S,T} \ | \ S \cup T \mbox{ is a partition of } [n] \}.$$
\end{prop}

\begin{ex}For $n =3,$ we have
\begin{eqnarray*}
\sC_{\sM} = &\{& P_{\{1\},\{2,3\}}=x_{1,2}x_{1,3}- x_{2,1}x_{3,1}, \\
& & P_{\{2,3\},\{1\}}=x_{2,1}x_{3,1} - x_{1,2}x_{1,3}, \\
& & P_{\{2\},\{1,3\}}=x_{2,1}x_{2,3} - x_{1,2}x_{3,2}, \\
& & P_{\{1,3\},\{2\}}=x_{1,2}x_{3,2} - x_{2,1}x_{2,3}, \\
& & P_{\{3\},\{1,2\}}=x_{3,1}x_{3,2} - x_{1,3}x_{2,3}, \\
& & P_{\{1,2\},\{3\}}=x_{1,3}x_{2,3} - x_{3,1}x_{3,2} \}.
\end{eqnarray*}
\end{ex}

We break the proof of Proposition \ref{circ} into several
lemmas. Before we state and prove the lemmas, we give a
formula for the entries in $$\u \sM = \sum_{i,j \in [n], i \neq j}
\u(i,j) \overline{M}_{i,j}.$$ For any $i,j \in [n-1],$ at most three
members of $\sM$ are nonzero at $(i,j)$-entry:
$\overline{M}_{i,j}(i,j) = 1$ (this one does not exist if $i = j$),
$\overline{M}_{i,n}(i,n) = -1,$ and $\overline{M}_{n,j}(n,j) = -1.$
Hence,
\begin{equation*}
(\u \sM)(i,j) = \begin{cases} - \u(i,n) - \u(n,j) & i = j; \\ \u(i,j)- \u(i,n) - \u(n,j) & i \neq j.
\end{cases}
\end{equation*}
Therefore, we have the following lemma.
\begin{lem}\label{uM}
\begin{equation*}
\u \sM = 0 \mbox{ if and only if } \begin{cases} \u(i,n) + \u(n,i) = 0, & \forall i \in [n-1] ;\\
\u(i,j) - \u(i,n) - \u(n,j) = 0,  & \forall i \neq j \in [n-1].
\end{cases}
\end{equation*}
\end{lem}

\begin{lem}\label{l1}
For any partition of $[n] = S \cup T,$ we have that
$\u_{S,T} \sM  = 0.$
Hence $P_{S,T}$ is in the toric ideal $I_{\sM}.$
\end{lem}
\begin{proof}
It directly follows from (\ref{anti}), (\ref{jacob}), and Lemma \ref{uM}.
\end{proof}

\begin{lem}\label{l2}
For any nonzero $\u \in \Z^{n(n-1)}$ satisfying $\u \sM = 0,$ i.e., $\x^{\u^+} - \x^{\u^-} \in I_{\sM},$ there exists a partition of $[n] = S \cup T,$ so that $\supp(\u_{S,T}) \subset \supp(\u).$
\end{lem}

\begin{proof}
We first show that there exists $t \in [n],$ such that either $(t,n)$
or $(n,t)$ is in the support $\supp(\u)$ of $\u.$ Let $(i,j) \in
\supp(\u),$ if either $i$ or $j$ is $n$, then we are done. Otherwise,
by Lemma \ref{uM}, we must have either $(i,n)$ or $(n,j)$ in
$\supp(\u).$

By Lemma \ref{uM} again, we conclude that $(t,n) \in \supp(\u)$ if and
only if $(n,t) \in \supp(\u).$ Let $T = \{t \ | \ (t, n) \in
\supp(\u)$ and/or $ (n, t) \in \supp(\u) \}$ and $S = [n] \setminus
T.$ Both $S$ and $T$ are nonempty. Thus $S \cup T$ is a partition of
$[n].$ We will show that $S \cup T$ is the partition needed to finish
the proof.

$\supp(\u_{S,T}) = \{ (s,t) \ | \ s \in S, t \in T\} \cup \{ (t,s) \ |
\ s \in S, t \in T\}.$ Hence, we need to show that $\forall s \in S,
\forall t \in T,$ both $(s,t)$ and $(t,s)$ are in $\supp(\u).$ If $s =
n,$ it follows immediately from the definition of $T.$ If $s \neq n,$
$(s,n) \not\in \supp(\u)$ since $s \not\in T.$ Therefore, $(\u
\sM)(s,t) = \u(s,t) - \u(n,t),$ which implies that $(s,t) \in
\supp(\u).$ We can similarly show that $(t,s) \in \supp(\u)$ as well.
\end{proof}

\begin{lem}\label{l3}
Let $\u \in \Z^{n(n-1)}$ satisfying $\u \sM = 0,$ and $\supp(\u) =
\supp(\u_{S,T})$ for some partition of $[n] = S \cup T,$ then $\exists
c \in \Z$ such that $\u = c \u_{S,T}.$
\end{lem}

\begin{proof}
Because $\u_{S,T} = - \u_{T,S},$ we can assume that $n \in S.$ Fix
$t_0 \in T,$ and let $c := \u(n, t_0),$ we will show that $\u = c
\u_{S,T}.$ Basically, we need to show that $\forall s \in S$ and
$\forall t \in T,$ $\u(s,t) = \u(n,t_0)$ and $\u(t,s) = - \u(n, t_0).$
We will show it case by case, by using Lemma \ref{uM} and the facts
$\u(s, n) =\u(n,s) = 0$ when $s \neq n$ and $\u(t_0, t) = 0$ when $t
\neq t_0.$
\begin{itemize}
\item If $s = n, t = t_0:$ $\u(s,t) = \u(n, t_0)$ and $\u(t,s) = \u(t_0, n) = - \u(n, t_0).$
\item If $s = n, t \neq t_0:$ $\u(s,t) = \u(n, t) = \u(t_0, t) - \u(t_0, n) = \u(n, t_0)$ and $\u(t,s) = \u(t,n) = -\u(n,t) = -\u(n, t_0).$
\item If $s \neq n, t = t_0:$ $\u(s,t) = \u(s, t_0) = \u(s,n) + \u(n,t_0) = \u(n, t_0)$ and $\u(t,s) = \u(t_0,s) = \u(t_0,n) + \u(n,s) = -\u(n, t_0).$
\item If $s \neq n, t \neq t_0:$ $\u(s,t) = \u(s,n) + \u(n,t) = \u(n, t_0)$ and $\u(t,s) = \u(t,n) + \u(n,s) = -\u(n, t_0).$
\end{itemize}
\end{proof}

\begin{proof}[Proof of Proposition \ref{circ}]
By Lemma \ref{l1}, Lemma \ref{l2} and Lemma \ref{l3}, we know that
$$\sC_{\sM} \subset \{ P_{S,T} \ | \ S \cup T \mbox{ is a partition of
} [n] \}.$$ Now we only need to show that for any partition $S \cup
T,$ there does not exist another partition $S' \cup T'$ such that
$\supp(\u_{S',T'})$ is strictly contained in $\supp(\u_{S,T}).$
Suppose we have such two partitions and let $(i,j) \in \supp(\u_{S,T})
\setminus \supp(\u_{S',T'}).$ Then $i$ and $j$ are both in $S'$ or
$T'.$ Without loss of generality, we assume they are both in $S'.$ Let
$t \in T',$ then $(i,t)$ and $(j,t)$ are both in the support of
$\u_{S',T'},$ thus in the support of $\u_{S,T}.$ But the fact that
$(i,j) \in \supp(S,T)$ indicates that one of $i$ and $j$ is in $S$ and
the other one is in $T.$ Wherever $t$ is in, we cannot have both
$(i,t)$ and $(j,t)$ are in the support of $\u_{S,T}.$ Therefore, we
proved that each $P_{S,T}$ is a circuit.
\end{proof}

\begin{cor}\label{gb}
For any $\ell \in [n],$ $$Gr_{\ell} := \{ P_{S,T} \ | \ S \cup T
\mbox{ is a partition of $[n]$ s.t. } \ell \in S \}$$ is a Gr\"obner
basis of $\sM$ with respect to any term order $<$ satisfying $x_{\ell,
j} > x_{i, k},$ for any $i \neq \ell.$ Thus, the set of initial
monomials of the elements in $Gr_{\ell}$ are
$$Ini(Gr_{\ell}) := \{ \prod_{s \in S, t \in T} x_{s,t} \ | \  S \cup T \mbox{ is a partition of $[n]$ s.t. } \ell \in S \}.$$
\end{cor}
\begin{ex}For $n =3, \ell = 3:$
\begin{eqnarray*}
Gr_{\ell} = &\{&
 P_{\{2,3\},\{1\}}=x_{2,1}x_{3,1} - x_{1,2}x_{1,3}, \\
& & P_{\{1,3\},\{2\}}=x_{1,2}x_{3,2} - x_{2,1}x_{2,3}, \\
& & P_{\{3\},\{1,2\}}=x_{3,1}x_{3,2} - x_{1,3}x_{2,3} \}
\end{eqnarray*}
and $$Ini(Gr_{\ell}) = \{x_{2,1}x_{3,1}, x_{1,2}x_{3,2},
x_{3,1}x_{3,2} \}.$$
\end{ex}

%
Recall that $\arb$ is the set of all $\ell$-arborescences on $[n].$
For any $T \in \arb$, we define the {\it support} of $T$ to be  $\supp(T) := \{ (i,j) \
| \ \mbox{ $i$ is the parent of $j$ in $T$}\}$, and let $\sM(T) = \{
\overline{M}_{i,j} \ | \ (i,j) \in \supp(T)\}$ be the corresponding
subset of $\sM$ defined in Lemma \ref{uni}. (Note the support of $T$ actually is the same as the edge set $E(T)$ of $T.$ We call it support here to be consistent with the definitions of other supports.)

\begin{prop} For any arborescence $T$ on $[n],$ we define $\overline{D}_T$ to be the cone generated by the rays in the set $\sM \setminus \sM(T),$ i.e., $D_T = \cone(\sM \setminus \sM(T)).$

Fix any $\ell \in [n],$ the term order and Gr\"obner basis described
in Corollary \ref{gb} give us a triangulation of $\overline{D}_n:$
$$Tri_{\ell}:= \{ \overline{D}_T \ | \ T \in \arb \}.$$
\end{prop}
\begin{proof}

From the theory of Gr\"obner bases of toric ideals in Theorem \ref{bernds}, the maximal simplices are given by the set
of transversals, all minimal sets $\sigma \subset \{(i,j) \ | \ i
\neq j \in [n] \}$ such that $\sigma \cap \supp(m) \neq \emptyset,
\forall m \in Ini(Gr_{\ell})$. Now due to the fact that each of the
initial monomials are in bijection to the cuts of on the complete
graph, the transversals are indeed given by all possible
arborescences
$$\{ \supp(T) \ | \ T \in \arb \} .$$ One direction is easy: given any
arborescence $T$ on $[n]$ with root $\ell$, one sees that $\supp(T)$
is a transversal. We show the other direction: if given a transversal
$\sigma,$ we can draw a directed graph $G_{\sigma}$ according to
$\sigma,$ i.e., $\supp(G_{\sigma}) = \sigma.$ We let $T$ be the set of
all $i$'s such that there does not exist a directed path from $\ell$
to $i.$ $T$ is empty, because otherwise $m = \prod_{s \not\in T, t \in
T} x_{s,t} \in Ini(Gr_{\ell})$ but $\sigma \cap \supp(m) = \emptyset.$
Therefore, for any vertex $i,$ there exists a directed path from
$\ell$ to $i.$ This implies that there is an $\ell$-arborescence as a
subgraph of $G_{\sigma}.$ However, by the minimality of $\sigma,$
$G_{\sigma}$ has to be this arborescence.

Finally, from Theorem \ref{bernds} we know
that the complement of these transversals are precisely the set of
simplices of the triangulation.
\end{proof}

\begin{ex} For $n=3, \ell =3,$
there are only three trees for $K_3$, thus the three
$3$-arborescences $T_A,T_B,T_C$ for $K_3$ are depicted in Figure
\ref{trees}.

\begin{figure}[h!]
    \begin{center}
        \includegraphics[scale=.40]{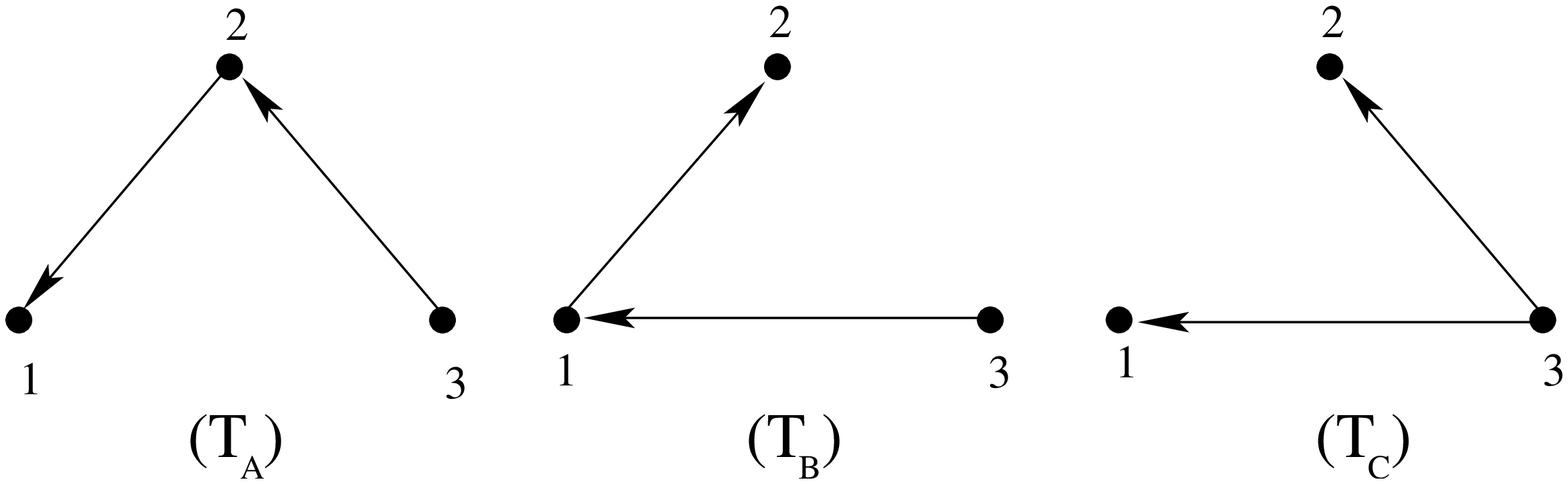}
    \end{center}
\caption{3-arborescences} \label{trees}
\end{figure}
\begin{eqnarray*}
Tri_{\ell} = &\{& \overline{D}_{T_A} = \cone(\sM \setminus \sM(T_A)) =
\cone(\{\overline{M}_{1,3},\overline{M}_{2,3},\overline{M}_{3,1},\overline{M}_{1,2}
\}) \\
& &\overline{D}_{T_B} = \cone(\sM \setminus \sM(T_B)) =
\cone(\{\overline{M}_{1,3},\overline{M}_{2,3},\overline{M}_{3,2},\overline{M}_{2,1}
\}) \\
& &\overline{D}_{T_C} = \cone(\sM \setminus \sM(T_C)) =
\cone(\{\overline{M}_{1,3},\overline{M}_{2,3},\overline{M}_{1,2},\overline{M}_{2,1}
\}) \},
\end{eqnarray*}
where $\overline{M}_{i,j}$ is defined as in Example \ref{exM1}.
\end{ex}

\subsection{Step 3: The dual cone to $\overline{D}_T$}

We have given triangulations $Tri_{\ell}$ of $\overline{D}_n.$ By Lemma \ref{uni}, we know this gives a decomposition of $\overline{D}_n$ into a set of unimodular cones $\overline{D}_T$, one for each arborescence $T$. Hence we can proceed to find the dual cone to each $\overline{D}_T$ inside $W_n$.

%



Recall $V_n$ is the subspace spanned by the supporting cone $C_n$ at the vertex $I$ and can be described by \eqref{sub}.
We will define cone $C_T,$ for each $T \in \arb,$ in the subspace $V_n,$ then show $\overline{C}_T := \phi(C_T)$ is the dual cone to $\overline{D}_T.$

\begin{defn}  For any 
directed edge $e = (s,t),$ ($s$ is pointed to $t$,) we define the {\bf
  weight} of $e$ to be the $n \times n$ matrix $w(e) \in \R^{n^2}$,
whose $(s,t)$-entry is $1,$ $(t,t)$-entry is $-1,$ and all the
remaining entries are zero.

Given $T$ an arborescence on $[n]$ with root $\ell,$ let $v$ be a
vertex of $T,$ then there is a unique path from $\ell$ to $v$, we
define the {\bf weight} $w^T(v)$ of $v$ with respect to $T$ to be the
summation of the weights of all the edges on this path.
\end{defn}

\begin{defn}\label{dualrays}
Let $T$ be an arborescence on $[n]$ with root $\ell$.
For each directed edge $e = (i,j)$ not in $T$, i.e., $e \not\in E(T),$ we define 
$$W^{T,e} := w^T(s) - w^T(t) + w(e).$$
More precisely, the entries of $W^{T,e}$ are
$$W^{T,e}(i,j) = \begin{cases}1, &  \mbox{if $i \not= j$ and $(i,j) \in cycle(T+e)$ has the same orientation as $e$},  \\
        -1, &  \mbox{if $i \not= j$ and $(i,j) \in cycle(T+e)$ has the opposite orientation as $e$}, \\
         -1, & \mbox{if $i=j$  and $i$ is a vertex in two edges of $cycle(T+e)$ with} \\
 & \mbox{both edges having same orientation as $e$.},    \\
         1, & \mbox{if $i=j$  and $i$ is a vertex in two edges of $cycle(T+e)$, with} \\
& \mbox{both edges having opposite orientation of $e$.},    \\
        0, &  \mbox{in all other cases}.
        \end{cases},$$
where $cycle(T+e)$ denote the unique cycle created by adding $e$ to $T$

Let $C_T$ be the cone generated by the set of rays $\{ W^{T,e} \ | \ e \not\in E(T) \}$ and $\overline{C}_T$ be its projection under $\phi$ (the map
that ignores the last row and last column of an $n \times n$ matrix):
$$C_T := \cone(\{ W^{T,e} \ | \ e \not\in E(T) \}), \ \ \overline{C}_T := \phi(C_T).$$
\end{defn}

\begin{prop}
\begin{enumerate}
\item Each $W^{T,e}$ is in the subspace $$V_n = \{ M \in \R^{n^2} \ | \ \sum_{k=1}^n M(i,k) = \sum_{k=1}^n M(k,j) = 0, \forall i, j\}.$$ Hence, $C_T$ is in $V_n.$ 
\item $\overline{C}_T$ is the dual cone to $\overline{D}_T$ in the vector space $W_n = \R^{(n-1)^2}.$
\end{enumerate}
\end{prop}

\begin{proof}\begin{enumerate}
\item We observe that for each row or column of $W^{T,e},$ there are either
one $1$, one $-1$ and  the other entries are zeros or all entries are zeros. 
\item $\overline{C}_T$ is the cone generated by the set of rays $\{ \phi({W}^{T,e}) \ | \ e \not\in E(T) \},$ and $\overline{D}_T$ is the cone generated by the set of rays $\{ \overline{M}_{i,j} \ | \ (i, j) \not\in \supp(T) \}.$ Recall $\phi$ is the map that ignores the entries in the last column and the last row of a matrix in $V_n \subset \R^{n^2}.$ Hence,  we have
$$\phi(W^{T,e})(k, \ell) = W^{T,e}(k, \ell), \forall 1 \le k, \ell \le n-1.$$
To check whether $\overline{C}_T$ is the dual cone to $\overline{D}_T$, it is enough to check for any directed edge $e = (s,t) \not\in E(T)$ and any $(i,j) \not\in \supp(T)$, we have $\langle \phi(W^{T,e}), \overline{M}_{i,j} \rangle$ is positive when $(i,j) = (s,t)$ and is $0$ otherwise. In fact, we
will show that $\langle \phi(W^{T,e}), \overline{M}_{i,j} \rangle = \delta_{(i,j),(s,t)}.$ 
There are three situations.
\begin{itemize}
\item If $i \neq n$ and $j \neq n,$ then $\langle \phi(W^{T,e}), \overline{M}_{i,j} \rangle  = \phi(W^{T,e})(i,j) = W^{T,e}(i,j).$
\item If $i = n$ and $j \neq n,$ then $\langle \phi(W^{T,e}), \overline{M}_{i,j} \rangle  = \sum_{k=1}^{n-1}(- \phi(W^{T,e})(k,j))= \sum_{k=1}^{n-1}(- W^{T,e}(k,j)) = W^{T,e}(n,j) = W^{T,e}(i,j).$
\item If $i \neq n$ and $j = n,$ similarly we have $\langle \phi(W^{T,e}), \overline{M}_{i,j} \rangle  = W^{T,e}(i,j).$
\end{itemize}
Hence, for every situation $\langle \phi(W^{T,e}), \overline{M}_{i,j} \rangle = W^{T,e}(i,j).$ However, since the only edge in $cycle(T+e)$ not in $T$ is $e,$ $W^{T,e}(i,j) = \delta_{(i,j),(s,t)}.$
\end{enumerate}
\end{proof}

\begin{ex}\label{exW} When $n = 3, \ell =3,$ as before we will 
present $W^{T,e}$ as a row vector,
which is just the first, second and last row of the matrix in order.
For the $3$-arborescence $T_A$ in Figure \ref{trees}, we have four
directed edges to be added, the edges $(1,2),(1,3),(3,1)$ and
$(2,3).$

$\begin{array}{cccccccccccc}
W^{T_A,(1,2)}: & -1 & 1 & 0 & &  1 & -1 & 0 & & 0 & 0 & 0     \\
W^{T_A,(1,3)}: & -1 & 0 & 1 & &  1 & -1 & 0 & & 0 & 1 & -1     \\
W^{T_A,(2,3)}: & 0 & 0 & 0 & &  0 & -1 & 1 & & 0 & 1 & -1  \\
W^{T_A,(3,1)}: & 0 & 0 & 0 & &  -1 & 1 & 0 & & 1 & -1 & 0
\end{array}.$

Similarly we have edges $(1,3), (2,1), (2,3)$ and $(3,2)$ to be
added onto the $3$-arborescence $T_B$ in Figure \ref{trees} and edges
$(1,2), (1,3), (2,1)$ and $(2,3)$ for the $3$-arborescence $T_C.$

$\begin{array}{cccccccccccc}
W^{T_B,(1,3)}: & -1 & 0 & 1 & & 0 & 0 & 0 & & 1 & 0 & -1 \\
W^{T_B,(2,1)}: & -1 & 1 & 0 & & 1 & -1 & 0 & & 0 & 0 & 0 \\
W^{T_B,(2,3)}: & -1 & 1 & 0 & & 0 & -1 & 1 & & 1 & 0 & -1 \\
W^{T_B,(3,2)}: & 1 & -1 & 0 & & 0 & 0 & 0 & & -1 & 1 & 0 \\
W^{T_C,(1,2)}: & -1 & 1 & 0 & & 0 & 0 & 0 & & 1 & -1 & 0 \\
W^{T_C,(1,3)}: & -1 & 0 & 1 & & 0 & 0 & 0 & & 1 & 0 & -1 \\
W^{T_C,(2,1)}: & 0 & 0 & 0 & & 1 & -1 & 0 & & -1 & 1 & 0 \\
W^{T_C,(2,3)}: & 0 & 0 & 0 & & 0 & -1 & 1 & & 0 & 1 & -1
\end{array}$
\end{ex}

\subsection{Step 4: The multivariate generating function of $C_n$}
Because each $\overline{D}_T$ in the triangulation of $\overline{D}_n$ is unimodular, so
is the dual cone $\overline{C}_T$ of $\overline{D}_T.$ By Lemma \ref{goodform}, we conclude that $C_T$ is unimodular as well and the following proposition:

\begin{prop}
Fixing $\ell \in [n],$ the multivariate generating function of
$C_n$ is given by
\begin{equation}\label{eqI}
f({C}_n,\z)= \sum_{T \in \arb} \prod_{e \notin E(T)}
\frac{1}{(1-\prod \z^{W^{T,e}})}.
\end{equation}
\end{prop}
One observes that Equation \eqref{eqI} is independent of the choice
of $\ell.$ Thus we have the following equality.
\begin{cor} \label{twoells}
For any $\ell_1, \ell_2 \in [n],$
$$\sum_{T \in {\bf Arb}(\ell_1, n)}  \prod_{e \notin E(T)}
\frac{1}{(1-\prod \z^{W^{T,e}})} = \sum_{T \in {\bf Arb}(\ell_2, n)}
\prod_{e \notin E(T)} \frac{1}{(1-\prod \z^{W^{T,e}})}.$$
\end{cor}

\section{A rational function formula for $f(tB_n,\z)$} \label{formulalatticepts}



In the last section, we obtained a formula for the multivariate
generating function of the supporting cone $C_n$ of the vertex $I$
of $B_n.$ Because of the symmetry of vertices of the
Birkhoff polytope we can get the MFG of the supporting cone of any other vertex of $B_n.$

\begin{cor} \label{eqsigma}
The multivariate generating function for the lattice points of
the supporting cone $C_n(\sigma)$ at the vertex $\sigma$, for $\sigma$ a permutation in $S_n$, is given by

\begin{equation} \label{vertexequ}
f(C_n(\sigma),\z)= \sum_{T \in \arb} \prod_{e \notin
E(T)} \frac{1}{(1-\prod \z^{W^{T,e}\sigma})},
\end{equation}
where $W^{T,e} \sigma$ is the matrix obtained from usual matrix multiplication of $W^{T,e}$ and the permutation matrix $\sigma.$
\end{cor}

\begin{proof}[Proof of Theorem \ref{main}] Note that for any positive integer $t,$ the supporting cone of $tB_n$ at vertex $t\sigma$ is still the same supporting cone $C_n(\sigma)$ of $B_n$ at the vertex $\sigma.$
Then the theorem follows from  Corollary \ref{intepoly} and Corollary \ref{eqsigma}.
\end{proof}

We conclude this section with an example of Theorem \ref{main} for
our running example.
\begin{ex}
When $n=3, \ell=3,$ the three $3$-arborescences are shown in Figure
\ref{trees}. In example \ref{exW}, we have already calculated
$W^{T,e}$'s. By plugging them in, we get the three parts of the
products of rational functions contributing to $f(C_3, \z):$

\begin{eqnarray*} \label{atreeequ}
\prod_{e \notin E(T_A)} \frac{1}{(1-\prod \z^{W^{T_A,e}})} &=&
\frac{1}{1-z_{1,2}z_{2,1}z_{1,1}^{-1}z_{2,2}^{-1}} \times
\frac{1}{1-z_{1,3}z_{3,2}z_{2,1}z_{1,1}^{-1}z_{2,2}^{-1}z_{3,3}^{-1}}
\times \\
& &\frac{1}{1-z_{2,3}z_{3,2}z_{2,2}^{-1}z_{3,3}^{-1}} \times
\frac{1}{1-z_{2,2}z_{3,1}z_{2,1}^{-1}z_{3,2}^{-1}}, \nonumber
\end{eqnarray*}
\begin{eqnarray*} \label{btreeequ}
\prod_{e \notin E(T_B)} \frac{1}{(1-\prod \z^{W^{T_B,e}})} &=&
\frac{1}{1-z_{1,3}z_{3,1}z_{1,1}^{-1}z_{3,3}^{-1}} \times
\frac{1}{1-z_{1,2}z_{2,1}z_{1,1}^{-1}z_{2,2}^{-1}} \times \\
&
&\frac{1}{1-z_{2,3}z_{3,1}z_{1,2}z_{1,1}^{-1}z_{2,2}^{-1}z_{3,3}^{-1}}
\times \frac{1}{1-z_{1,1}z_{3,2}z_{1,2}^{-1}z_{3,1}^{-1}}, \nonumber
\end{eqnarray*} and
\begin{eqnarray*} \label{ctreeequ}
\prod_{e \notin E(T_C)} \frac{1}{(1-\prod \z^{W^{T_C,e}})} &=&
\frac{1}{1-z_{1,2}z_{3,1}z_{1,1}^{-1}z_{3,2}^{-1}} \times
\frac{1}{1-z_{1,3}z_{3,1}z_{1,1}^{-1}z_{3,3}^{-1}} \times\\
& &\frac{1}{1-z_{2,1}z_{3,2}z_{2,2}^{-1}z_{3,1}^{-1}} \times
\frac{1}{1-z_{2,3}z_{3,2}z_{3,3}^{-1}z_{2,2}^{-1}}. \nonumber
\end{eqnarray*}
Thus, $\z^{tI} f(C_3, \z)$ equals the sum of
the three rational functions multiplied by $(z_{1,1}^tz_{2,2}^tz_{3,3}^t).$

In order to compute the same for other vertices we simply permute
the results:
\begin{eqnarray*}
& &\z^{t\sigma} f(C_3(\sigma), \z) = (z_{1, \sigma(1)}^t z_{2, \sigma(2)}^t
z_{3, \sigma(3)}^t) \times \\
& &(
\frac{1}{1-z_{1,\sigma(2)}z_{2,\sigma(1)}z_{1,\sigma(1)}^{-1}z_{2,\sigma(2)}^{-1}}
\times
\frac{1}{1-z_{1,\sigma(3)}z_{3,\sigma(2)}z_{2,\sigma(1)}z_{1,\sigma(1)}^{-1}z_{2,\sigma(2)}^{-1}z_{3,\sigma(3)}^{-1}}
\times \\
&
&\frac{1}{1-z_{2,\sigma(3)}z_{3,\sigma(2)}z_{2,\sigma(2)}^{-1}z_{3,\sigma(3)}^{-1}}
\times
\frac{1}{1-z_{2,\sigma(2)}z_{3,\sigma(1)}z_{2,\sigma(1)}^{-1}z_{3,\sigma(2)}^{-1}}
+ \\
&
&\frac{1}{1-z_{1,\sigma(3)}z_{3,\sigma(1)}z_{1,\sigma(1)}^{-1}z_{3,\sigma(3)}^{-1}}
\times
\frac{1}{1-z_{1,\sigma(2)}z_{2,\sigma(1)}z_{1,\sigma(1)}^{-1}z_{2,\sigma(2)}^{-1}} \times \\
&
&\frac{1}{1-z_{2,\sigma(3)}z_{3,\sigma(1)}z_{1,\sigma(2)}z_{1,\sigma(1)}^{-1}z_{2,\sigma(2)}^{-1}z_{3,\sigma(3)}^{-1}}
\times
\frac{1}{1-z_{1,\sigma(1)}z_{3,\sigma(2)}z_{1,\sigma(2)}^{-1}z_{3,\sigma(1)}^{-1}}
+ \\
&
&\frac{1}{1-z_{1,\sigma(2)}z_{3, \sigma(1)}z_{1,\sigma(1)}^{-1}z_{3,\sigma(2)}^{-1}}
\times
\frac{1}{1-z_{1,\sigma(3)}z_{3,\sigma(1)}z_{1,\sigma(1)}^{-1}z_{3,\sigma(3)}^{-1}} \times\\
&
&\frac{1}{1-z_{2,\sigma(1)}z_{3,\sigma(2)}z_{2,\sigma(2)}^{-1}z_{3,\sigma(1)}^{-1}}
\times
\frac{1}{1-z_{2,\sigma(3)}z_{3,\sigma(2)}z_{3,\sigma(3)}^{-1}z_{2,\sigma(2)}^{-1}}).
\end{eqnarray*}

Finally, the summation of all six $\z^{t\sigma} f(C_3(\sigma), \z)$ gives
$f(tB_3, \z).$
\end{ex}

\section{The Coefficients of the Ehrhart polynomial of the Birkhoff polytope}
\label{coeffsehrhart}
In section $5.2$ of \cite{BarviPom}, Barvinok and Pommersheim derive
a formula for the number of lattice points of a given integral convex
polytope $P$ in terms of {\it Todd polynomial} by residue
computation of the MGF of $P.$ When $P$ is an integral polytope,
their formula explicitly indicates formulas for the coefficients
of the Ehrhart polynomial $e(P,t)$ of $P.$ Especially, this gives us
a formula for the volume $\vol(P)$ of $P,$ applying it we can get
Theorem \ref{volume}. We start this section by briefly 
recalling related results in \cite{BarviPom}.

\begin{defn} \label{todddef}
Consider the function
$$G(\tau; \xi_1, \ldots, \xi_d)=\prod_{i=1}^d
{\tau \xi_i \over 1-\exp(-\tau \xi_i)}$$ in $d+1$ (complex)
variables $\tau$ and $\xi_1, \ldots, \xi_l$. The function $G$ is
analytic in a neighborhood of the origin $\tau=\xi_1= \ldots
=\xi_d=0$ and therefore there exists an expansion
$$G(\tau; \xi_1, \ldots, \xi_d)=\sum_{j=0}^{+\infty} \tau^j
\td_j(\{\xi_i| 1 \le i \le d\}),$$ where $\td_j(\{\xi_i | 1\le i \le
d\})=\td_j(\xi_1,\xi_2,\dots,\xi_d)$ is a homogeneous polynomial of
degree $j$, called the $j$-th {\bf Todd polynomial} in $\xi_1, \ldots,
\xi_d$. It is well-know that $\td_j(\{\xi_i | 1 \le i \le d\})$ is a
symmetric polynomial with rational coefficients. See page 
110 in \cite{Fulton} for more information on Todd polynomials.
\end{defn}

\begin{ex} Here are the first three Todd polynomials when $d=3$:
$$
\td_3(x_1,x_2,x_3)=(1/24)\, \left( {x_1}+{x_2}+{x_3} \right)  \left( {x_1
}{x_2}+{x_2}{x_3}+{x_3}{x_1} \right),$$
$$\td_2(x_1,x_2,x_3)=(1/12){{x_2}}^{2}+(1/4){x_3}{x_1}+(1/12){{x_3}}^{2
}+(1/12){{x_1}}^{2}+(1/4){x_2}{x_3}+(1/4){x_1}{x_2},$$
$$\td_1(x_1,x_2,x_3)=
(1/2)\,{x_1}+(1/2)\,{x_2}+(1/2)\,{ x_3},\ \mbox{and as usual} \quad \td_0(x_1,x_2,x_3)=1.$$
\end{ex}

\begin{lem}\label{latticecount}(See Algorithm 5.2 in \cite{BarviPom}.)
Suppose $P \subset \R^N$ is $d$-dimensional integral polytope and the
multivariate generating function of $P$ is given by
\begin{equation}\label{mgf}
f(P, \z) = \sum_{i} \epsilon_i
\frac{\z^{a_i}_i}{(1-\z^{b_{i,1}})\cdots(1-\z^{b_{i,d}})},
\end{equation}
 where
$\epsilon_i = \{-1, 1\},$ $a_i,b_{i,1}, \dots, b_{i,d} \in \Z^N$, 
the $a_i$'s are all vertices (with multiple occurrences) of $P.$
and $\cone(b_{i,1}, \dots, b_{i,d})$ is unimodular, for each $i.$
For any choice of $c \in \R^N$ such that $\langle c, b_{i,j} \rangle \neq 0$ for each $i$ and $j,$ we have a formula for the number of lattice points in $P:$
\begin{equation}\label{lp}
|P \cap \Z^N| = \sum_i \frac{\epsilon_i}{\prod_{j=1}^d \langle c ,
b_{i,j} \rangle} \sum_{k=0}^d \frac{(\langle c, a_i \rangle)^k}{k!}
\td_{d-k}(\langle c , b_{i,1} \rangle, \dots, \langle c , b_{i,d}
\rangle).
\end{equation}
\end{lem}

Indeed, if we make the substitution $x_i=\exp(\tau c_i)$ Formula
(\ref{mgf}) can be rewritten as
\begin{equation}
f(P, \z) = \frac{1}{\tau^d} \sum_{i} \epsilon_i
\frac{\tau^d \exp(\langle c, a_i \rangle}{(1-\exp(\langle c, b_{i,1} \rangle)\cdots(1-\exp(\langle c, b_{i,d} \rangle)}.
\end{equation}
Each fraction is a holomorphic function in a neighborhood of $\tau$
and the $d$-th coefficient of its Taylor series is a linear
combination of Todd polynomials. Thus its $d$-coefficient of the
Taylor series is
\begin{equation}
\frac{1}{(\langle c, b_{i,1} \rangle)\cdots (\langle c, b_{i,d} \rangle)} \sum_{k=0}^d \frac{ (\langle c, a_i \rangle)^k}{k!} \td_{d-k}(\langle c, b_{i,1} \rangle,\dots, \langle c, b_{i,d} \rangle).
\end{equation}
Formula (\ref{lp}) is the result of adding these contributions for each rational fraction summand.

It is clear that if Formula \eqref{mgf} is the
MGF of an integral polytope $P,$ then we have the MGF of any of 
its  dilations:
\begin{equation}\label{mgfdil}
f(tP, \z) = \sum_{i} \epsilon_i
\frac{\z^{ta_i}}{(1-\z^{b_{i,1}})\cdots(1-\z^{b_{i,d}})}.
\end{equation}
Hence, by using the Lemma \ref{latticecount}, we get the Ehrhart
polynomial of $P.$

\begin{lem}\label{vollem}
Suppose $P \subset \R^N$ is $d$-dimensional integral polytope and
the multivariate generating function of $P$ (produced by Barvinok's
algorithm) is given by \eqref{mgf}. For any choice of $c \in \R^N$
such that $\langle c, b_{i,j} \rangle \neq 0$ for each $i$ and $j,$
the Ehrhart polynomial of $P$ is
\begin{equation}\label{ehrhart}
e(P,t) = \sum_{k=0}^d \frac{t^k}{k!}\sum_i
\frac{\epsilon_i}{\prod_{j=1}^d \langle c , b_{i,j} \rangle}
{(\langle c, a_i \rangle)^k} \td_{d-k}(\langle c , b_{i,1} \rangle,
\dots, \langle c , b_{i,d} \rangle).
\end{equation}
In particular, we get a formula for the volume of $P:$
\begin{equation}\label{volform}
\vol(P) = \frac{1}{d!}\sum_i \epsilon_i \frac{(\langle c, a_i
\rangle)^d}{\prod_{j=1}^d \langle c , b_{i,j} \rangle} .
\end{equation}
\end{lem}

\begin{proof}
Formula \eqref{ehrhart} directly follows Lemma \ref{latticecount}
and our earlier discussion. Formula \eqref{volform} follows the
facts that the leading coefficient of $e(P,t)$ is $\vol(P)$ and the
$0$-th Todd polynomial is always the constant $1$, which can be shown from 
the Taylor expression of the function $G(\tau,\xi_1,\dots,\xi_d)$ defining
the Todd polynomials.
\end{proof}

\begin{proof}[Proof of Corollary \ref{volume}]
It follows Lemma \ref{vollem} and Theorem \ref{main}.
\end{proof}


To help our readers we wrote an interactive {\tt MAPLE} implementation
of Formula \eqref{volume} in the case of $B_3$. It is available at
\url{http://www.math.ucdavis.edu/~deloera/RECENT_WORK/volBirkhoff3}

Clearly, it would be desirable to apply a suitable variable
substitution of $c_{i,j}$ so that the expression of the volume has as
few terms as possible (preferably keeping size of $c_{i,j}$ small),
with the hope of speeding up the calculations or even in hope of
finding a purely combinatorial summation.  We leave this challenge to
the reader and conclude with a variable exchange that gives the volume
in just two variables (it is possible to leave it as a univariate
rational function from the substitution $c_{i,j}=it^j$). If we set
$c_{i,j}=s^it^j$ clearly there will be no cancellations.  For example
for the case $n=3$, the volume of $B_3$ equals.

{\tiny 
\noindent $
1/24\,{\frac { \left( st+{s}^{2}{t}^{2}+{s}^{3}{t}^{3} \right) ^{4}}{
 \left( s{t}^{2}+{s}^{2}t-st-{s}^{2}{t}^{2} \right)  \left( {s}^{2}{t}
^{3}+{s}^{3}{t}^{2}-{s}^{2}{t}^{2}-{s}^{3}{t}^{3} \right)  \left( s{t}
^{3}+{s}^{3}{t}^{2}+{s}^{2}t-st-{s}^{2}{t}^{2}-{s}^{3}{t}^{3} \right) 
 \left( {s}^{2}{t}^{2}+{s}^{3}t-{s}^{2}t-{s}^{3}{t}^{2} \right) }}+ \\
 1/24\,{\frac { \left( st+{s}^{2}{t}^{3}+{s}^{3}{t}^{2} \right) ^{4}}{
 \left( s{t}^{3}+{s}^{2}t-st-{s}^{2}{t}^{3} \right)  \left( {s}^{2}{t}
^{2}+{s}^{3}{t}^{3}-{s}^{2}{t}^{3}-{s}^{3}{t}^{2} \right)  \left( s{t}
^{2}+{s}^{3}{t}^{3}+{s}^{2}t-st-{s}^{2}{t}^{3}-{s}^{3}{t}^{2} \right) 
 \left( {s}^{2}{t}^{3}+{s}^{3}t-{s}^{2}t-{s}^{3}{t}^{3} \right) }}+ \\
 1/24\,{\frac { \left( s{t}^{2}+{s}^{2}t+{s}^{3}{t}^{3} \right) ^{4}}{
 \left( st+{s}^{2}{t}^{2}-s{t}^{2}-{s}^{2}t \right)  \left( {s}^{2}{t}
^{3}+{s}^{3}t-{s}^{2}t-{s}^{3}{t}^{3} \right)  \left( s{t}^{3}+{s}^{3}
t+{s}^{2}{t}^{2}-s{t}^{2}-{s}^{2}t-{s}^{3}{t}^{3} \right)  \left( {s}^
{2}t+{s}^{3}{t}^{2}-{s}^{2}{t}^{2}-{s}^{3}t \right) }}+\\
1/24\,{\frac {
 \left( s{t}^{2}+{s}^{2}{t}^{3}+{s}^{3}t \right) ^{4}}{ \left( s{t}^{3
}+{s}^{2}{t}^{2}-s{t}^{2}-{s}^{2}{t}^{3} \right)  \left( {s}^{2}t+{s}^
{3}{t}^{3}-{s}^{2}{t}^{3}-{s}^{3}t \right)  \left( st+{s}^{3}{t}^{3}+{
s}^{2}{t}^{2}-s{t}^{2}-{s}^{2}{t}^{3}-{s}^{3}t \right)  \left( {s}^{2}
{t}^{3}+{s}^{3}{t}^{2}-{s}^{2}{t}^{2}-{s}^{3}{t}^{3} \right) }}+\\
1/24\,{\frac { \left( s{t}^{3}+{s}^{2}t+{s}^{3}{t}^{2} \right) ^{4}}{
 \left( st+{s}^{2}{t}^{3}-s{t}^{3}-{s}^{2}t \right)  \left( {s}^{2}{t}
^{2}+{s}^{3}t-{s}^{2}t-{s}^{3}{t}^{2} \right)  \left( s{t}^{2}+{s}^{3}
t+{s}^{2}{t}^{3}-s{t}^{3}-{s}^{2}t-{s}^{3}{t}^{2} \right)  \left( {s}^
{2}t+{s}^{3}{t}^{3}-{s}^{2}{t}^{3}-{s}^{3}t \right) }}+\\1/24\,{\frac {
 \left( s{t}^{3}+{s}^{2}{t}^{2}+{s}^{3}t \right) ^{4}}{ \left( s{t}^{2
}+{s}^{2}{t}^{3}-s{t}^{3}-{s}^{2}{t}^{2} \right)  \left( {s}^{2}t+{s}^
{3}{t}^{2}-{s}^{2}{t}^{2}-{s}^{3}t \right)  \left( st+{s}^{3}{t}^{2}+{
s}^{2}{t}^{3}-s{t}^{3}-{s}^{2}{t}^{2}-{s}^{3}t \right)  \left( {s}^{2}
{t}^{2}+{s}^{3}{t}^{3}-{s}^{2}{t}^{3}-{s}^{3}{t}^{2} \right) }}+\\
1/24\,{\frac { \left( st+{s}^{2}{t}^{2}+{s}^{3}{t}^{3} \right) ^{4}}{
 \left( s{t}^{2}+{s}^{2}t-st-{s}^{2}{t}^{2} \right)  \left( s{t}^{3}+{
s}^{3}t-st-{s}^{3}{t}^{3} \right)  \left( {s}^{2}{t}^{3}+{s}^{3}t+s{t}
^{2}-st-{s}^{2}{t}^{2}-{s}^{3}{t}^{3} \right)  \left( st+{s}^{3}{t}^{2
}-s{t}^{2}-{s}^{3}t \right) }}+\\
1/24\,{\frac { \left( st+{s}^{2}{t}^{3}
+{s}^{3}{t}^{2} \right) ^{4}}{ \left( s{t}^{3}+{s}^{2}t-st-{s}^{2}{t}^
{3} \right)  \left( s{t}^{2}+{s}^{3}t-st-{s}^{3}{t}^{2} \right) 
 \left( {s}^{2}{t}^{2}+{s}^{3}t+s{t}^{3}-st-{s}^{2}{t}^{3}-{s}^{3}{t}^
{2} \right)  \left( st+{s}^{3}{t}^{3}-s{t}^{3}-{s}^{3}t \right) }}+\\
1/24\,{\frac { \left( s{t}^{2}+{s}^{2}t+{s}^{3}{t}^{3} \right) ^{4}}{
 \left( st+{s}^{2}{t}^{2}-s{t}^{2}-{s}^{2}t \right)  \left( s{t}^{3}+{
s}^{3}{t}^{2}-s{t}^{2}-{s}^{3}{t}^{3} \right)  \left( {s}^{2}{t}^{3}+{
s}^{3}{t}^{2}+st-s{t}^{2}-{s}^{2}t-{s}^{3}{t}^{3} \right)  \left( s{t}
^{2}+{s}^{3}t-st-{s}^{3}{t}^{2} \right) }}+\\
1/24\,{\frac { \left( s{t}^
{2}+{s}^{2}{t}^{3}+{s}^{3}t \right) ^{4}}{ \left( s{t}^{3}+{s}^{2}{t}^
{2}-s{t}^{2}-{s}^{2}{t}^{3} \right)  \left( st+{s}^{3}{t}^{2}-s{t}^{2}
-{s}^{3}t \right)  \left( {s}^{2}t+{s}^{3}{t}^{2}+s{t}^{3}-s{t}^{2}-{s
}^{2}{t}^{3}-{s}^{3}t \right)  \left( s{t}^{2}+{s}^{3}{t}^{3}-s{t}^{3}
-{s}^{3}{t}^{2} \right) }}+\\
1/24\,{\frac { \left( s{t}^{3}+{s}^{2}t+{s}
^{3}{t}^{2} \right) ^{4}}{ \left( st+{s}^{2}{t}^{3}-s{t}^{3}-{s}^{2}t
 \right)  \left( s{t}^{2}+{s}^{3}{t}^{3}-s{t}^{3}-{s}^{3}{t}^{2}
 \right)  \left( {s}^{2}{t}^{2}+{s}^{3}{t}^{3}+st-s{t}^{3}-{s}^{2}t-{s
}^{3}{t}^{2} \right)  \left( s{t}^{3}+{s}^{3}t-st-{s}^{3}{t}^{3}
 \right) }}+\\
 1/24\,{\frac { \left( s{t}^{3}+{s}^{2}{t}^{2}+{s}^{3}t
 \right) ^{4}}{ \left( s{t}^{2}+{s}^{2}{t}^{3}-s{t}^{3}-{s}^{2}{t}^{2}
 \right)  \left( st+{s}^{3}{t}^{3}-s{t}^{3}-{s}^{3}t \right)  \left( {
s}^{2}t+{s}^{3}{t}^{3}+s{t}^{2}-s{t}^{3}-{s}^{2}{t}^{2}-{s}^{3}t
 \right)  \left( s{t}^{3}+{s}^{3}{t}^{2}-s{t}^{2}-{s}^{3}{t}^{3}
 \right) }}+\\
 1/24\,{\frac { \left( st+{s}^{2}{t}^{2}+{s}^{3}{t}^{3}
 \right) ^{4}}{ \left( s{t}^{3}+{s}^{3}t-st-{s}^{3}{t}^{3} \right) 
 \left( {s}^{2}{t}^{3}+{s}^{3}{t}^{2}-{s}^{2}{t}^{2}-{s}^{3}{t}^{3}
 \right)  \left( s{t}^{2}+{s}^{3}t-st-{s}^{3}{t}^{2} \right)  \left( {
  s}^{2}t+{s}^{3}{t}^{2}-{s}^{2}{t}^{2}-{s}^{3}t \right) }}+\\
 1/24\,{
\frac { \left( st+{s}^{2}{t}^{3}+{s}^{3}{t}^{2} \right) ^{4}}{ \left( 
s{t}^{2}+{s}^{3}t-st-{s}^{3}{t}^{2} \right)  \left( {s}^{2}{t}^{2}+{s}
^{3}{t}^{3}-{s}^{2}{t}^{3}-{s}^{3}{t}^{2} \right)  \left( s{t}^{3}+{s}
^{3}t-st-{s}^{3}{t}^{3} \right)  \left( {s}^{2}t+{s}^{3}{t}^{3}-{s}^{2
}{t}^{3}-{s}^{3}t \right) }}+\\
1/24\,{\frac { \left( s{t}^{2}+{s}^{2}t+{
s}^{3}{t}^{3} \right) ^{4}}{ \left( s{t}^{3}+{s}^{3}{t}^{2}-s{t}^{2}-{
s}^{3}{t}^{3} \right)  \left( {s}^{2}{t}^{3}+{s}^{3}t-{s}^{2}t-{s}^{3}
{t}^{3} \right)  \left( st+{s}^{3}{t}^{2}-s{t}^{2}-{s}^{3}t \right) 
 \left( {s}^{2}{t}^{2}+{s}^{3}t-{s}^{2}t-{s}^{3}{t}^{2} \right) }}+\\
 1/24\,{\frac { \left( s{t}^{2}+{s}^{2}{t}^{3}+{s}^{3}t \right) ^{4}}{
 \left( st+{s}^{3}{t}^{2}-s{t}^{2}-{s}^{3}t \right)  \left( {s}^{2}t+{
s}^{3}{t}^{3}-{s}^{2}{t}^{3}-{s}^{3}t \right)  \left( s{t}^{3}+{s}^{3}
{t}^{2}-s{t}^{2}-{s}^{3}{t}^{3} \right)  \left( {s}^{2}{t}^{2}+{s}^{3}
{t}^{3}-{s}^{2}{t}^{3}-{s}^{3}{t}^{2} \right) }}+\\
1/24\,{\frac {
 \left( s{t}^{3}+{s}^{2}t+{s}^{3}{t}^{2} \right) ^{4}}{ \left( s{t}^{2
}+{s}^{3}{t}^{3}-s{t}^{3}-{s}^{3}{t}^{2} \right)  \left( {s}^{2}{t}^{2
}+{s}^{3}t-{s}^{2}t-{s}^{3}{t}^{2} \right)  \left( st+{s}^{3}{t}^{3}-s
{t}^{3}-{s}^{3}t \right)  \left( {s}^{2}{t}^{3}+{s}^{3}t-{s}^{2}t-{s}^
{3}{t}^{3} \right) }}+\\
1/24\,{\frac { \left( s{t}^{3}+{s}^{2}{t}^{2}+{s
}^{3}t \right) ^{4}}{ \left( st+{s}^{3}{t}^{3}-s{t}^{3}-{s}^{3}t
 \right)  \left( {s}^{2}t+{s}^{3}{t}^{2}-{s}^{2}{t}^{2}-{s}^{3}t
 \right)  \left( s{t}^{2}+{s}^{3}{t}^{3}-s{t}^{3}-{s}^{3}{t}^{2}
 \right)  \left( {s}^{2}{t}^{3}+{s}^{3}{t}^{2}-{s}^{2}{t}^{2}-{s}^{3}{
t}^{3} \right) }}. $
}

\section{Integration of polynomials and volumes of faces of the Birkhoff polytope} \label{facesofBn}

In this final section we look at two more applications of Theorem \ref{main}, going beyond the computation
of Ehrhart coefficients. 

The first application is to the integration of polynomials over $B_n$. The main observation is that, once we know a unimodular cone decomposition
for the supporting cones at all vertices of $B_n$, a formula for the integral (see Formula (\ref{intbrion})) follows from Brion's theorem on polyhedra \cite{barvivolume, Brion88}.

\begin{thm} 
Suppose $P \subset \R^N$ is a $d$-dimensional integral polytope and the
multivariate generating function of $P$ is given precisely by Formula \eqref{mgf}, i.e., we have full
knowledge of a unimodular cone decomposition for each of $P$'s supporting cones and its rays $b_{i,j}$ and vertices $a_i$. Then for any choice of $y \in \R^N$
such that $\langle y, b_{i,j} \rangle \neq 0$ for each $i$ and $j,$ we get a formula for the integral the $p$-th power of a linear form over $P$
\begin{equation}\label{intbrion}
\int_{P} {\langle y,x \rangle}^p dx = \frac{(-1)^d}{(p+1)(p+2)\dots(p+d)}  \sum_i \epsilon_i   \frac{(\langle y, a_i \rangle)^{p+d}}{\prod_{j=1}^d \langle y , b_{i,j} \rangle} .
\end{equation}
\end{thm}

Notice that although each term in the sum has poles, the poles cancel and the sum is 
an analytic function of $y$.
Since the $p$th-powers of linear forms generate the whole vector spaces of polynomials one obtains, from
Theorem \ref{main}, a formula of integration for the polynomial functions over $B_n$ (or for that matter,
for any integral polytope for which we understand their cone decomposition).

The next  application to the computation of Ehrhart polynomials of faces of $B_n$.
One can easily obtain from Theorem \ref{main} similar formulas for the nonnegative integral
semi-magic squares with structural zeros or forbidden entries
(i.e. fixed entries are equal zero). Note that any face
$F$ of $B_n$, being the intersection of finitely many facets, is determine uniquely by the set of entries forced to
take the value zero.  To obtain a generating function for the
dilations of a face $F$ of $B_n$, $f(tF,\z)$, we start from our formula
for $f(tB_n,z)$ in Theorem \ref{main}.  For those variables $x_{ij}$
mandated to be zero, we select a vector $\lambda$ with entries
$\lambda_{ij} \geq 0$ so that the substitution $x_{ij}:=s^{\lambda_{ij}}$ does not create a singularity (this
$\lambda$ exists e.g., by taking random values from the positive
orthant). Call $g(tF,\z,s)$ the result of doing this substitution on
$f(tB_n,\z)$.  We will eventually set $s=0$, but first, let us check
that will not give any singularities. On the numerator $s$ can only
appear with a nonnegative exponent. It can potentially create a
singularity if it appears in a factor of the denominator with negative
exponent. But, if this occurs, $s$ can be factored out and put with a
positive exponent in the numerator. Thus, we can safely resolve the
singularity.  Now, we set $s=0$ in $g(tF,\z,s)$, those terms that had a
power of $s$ in the numerator disappear. We obtain a multivariate sum
of rational functions that gives us only the desired lattice points
inside $tF$.  We have now two examples of this method. First we apply it to obtain a table
with the Ehrhart polynomials for  a (any) facet of $B_3, B_4, B_5, B_6$.

{\tiny
\[
\begin{array}{|c|l|}\hline
n & \mbox{Ehrhart polynomial of a facet of $B_n$}\\\hline
& \\
3 & 1 + \frac{11}{6} t + t^2 + \frac{1}{6}t^3\\
& \\\hline
& \\
4 & 1 + \frac{471}{140}t + \frac{1594}{315}t^2 + \frac{73}{16}t^3 +
\frac{161}{60}t^4 + \frac{83}{80}t^5 + \frac{61}{240}t^6 + \frac{1}{28}t^7
+ \frac{11}{5040}t^8\\
& \\\hline
& \\
5 & 1 + \frac{1752847}{360360}t + \frac{904325}{77616}t^2 +
\frac{147579347}{8072064}t^3 + \frac{8635681}{415800}t^4 +
\frac{6412937357}{359251200}t^5 + \frac{18455639}{1555200}t^6 +  \\
& \\
& \frac{1611167963}{261273600}t^7 + \frac{95702009}{38102400}t^8 + \frac{365214839}{457228800}t^9 + \frac{5561}{28350}t^{10} +
 \frac{52388227}{1437004800}t^{11} + \frac{42397}{8553600}t^{12} + \\
& \\
& \frac{4342517}{9340531200}t^{13} + \frac{22531}{838252800}t^{14} +
\frac{188723}{261534873600}t^{15}\\
& \\
\hline
& \\
6 & 1 + \frac{87450005}{13728792}t+ \frac{102959133218657}{4947307485120}t^2+ \frac{14843359499161}{322075353600}t^3+ \frac{230620928072832499}{3011404556160000}t^4+ \frac{237290485580450429}{2365321396838400}t^5+\\
& \\
&  \frac{15435462135033037}{144815595724800}t^6+ \frac{108878694347719}{1164067946496}t^7+ \frac{439368248888657369}{6402373705728000}t^8+ \frac{44766681773591807}{1054508610355200}t^9+\\
& \\
&  \frac{434798171323757}{19527937228800}t^{10}+ \frac{1047553900202141}{105450861035520}t^{11}+ \frac{250284934507924171}{66283398365184000}t^{12}+ \frac{28330897394929}{23176013414400}t^{13}+ \\
& \\
&  \frac{2229552439625171}{6628339836518400}t^{14}+ \frac{6610306048279}{84360688828416}t^{15}+ \frac{215934508972451}{14060114804736000}t^{16}+ \frac{2045239925737}{814847562547200}t^{17}+ \\
& \\
&  \frac{1729908621731}{5121898964582400}t^{18}+ \frac{21042914689}{572447531335680}t^{19}+\frac{6138921521069}{1946321606541312000}t^{20}+ \frac{139856666897}{681212562289459200}t^{21}+ \\
& \\
&  \frac{47580345877}{4995558790122700800}t^{22}+ \frac{4394656999}{15667888932657561600}t^{23}+ \frac{9700106723}{2462096832274759680000}t^{24}\\
& \\
\hline

\end{array}
\]
}

With the same method we have also computed, for the first time, the Ehrhart polynomials for
the Chan-Robbins-Yen Polytope $CRY_3$,$CRY_4$,$CRY_5$,$CRY_6$, and
$CRY_7$ (see \cite{chanrobbinsyuen}):

{\tiny
\[
\begin{array}{|c|l|}\hline
n & \mbox{Ehrhart polynomial of the Chan-Robbins-Yen polytopes}\\\hline
& \\
3 & 1 + \frac{11}{6} t + t^2 + \frac{1}{6}t^3\\
& \\\hline
& \\
4 & 1 + \frac{157}{60}t + \frac{949}{360}t^2 + \frac{4}{3}t^3 +
\frac{13}{36}t^4 + \frac{1}{20}t^5 + \frac{1}{360}t^6 \\
& \\ \hline
& \\
5 & 1 +\frac{2843}{840}t+ \frac{1087}{224}t^2 + \frac{16951}{4320}t^3 +
\frac{723869}{362880}t^4 + \frac{1927}{2880}t^5 + \frac{2599}{17280}t^6 +
\frac{113}{5040}t^7 + \frac{257}{120960}t^8 + \frac{1}{8640}t^9 +
\frac{1}{362880}t^{10} \\
& \\ \hline
& \\
6 & 1 + \frac{1494803}{360360}t+ \frac{15027247}{1965600}t^2+
\frac{361525133}{43243200}t^3+ \frac{364801681}{59875200}t^4+
\frac{45175393}{14370048}t^5+ \frac{4314659}{3628800}t^6+
\frac{4392257}{13063680}t^7+\\
& \\
& \frac{781271}{10886400}t^8+ \frac{75619}{6531840}t^9+
\frac{15257}{10886400}t^{10}+ \frac{22483}{179625600}t^{11}+
\frac{29}{3628800}t^{12}+ \frac{23}{66718080}t^{13} + \frac{1}{111196800}t^{14}
+ \\& \\ & \frac{1}{9340531200}t^{15}\\
& \\ \hline
& \\
7 & 1 +\frac{571574671}{116396280}t+ \frac{41425488163}{3760495200}t^2+
\frac{88462713645601}{5866372512000}t^3+
\frac{26256060764993}{1852538688000}t^4+
\frac{433329666631051}{44460928512000}t^5+ \\
& \\
&\frac{615428916451}{120708403200}t^6+
\frac{97984316095277}{47076277248000}t^7+
\frac{7939938012827}{11769069312000}t^8+
\frac{66150911695291}{376610217984000}t^9+ \\
& \\
& \frac{71471423831}{1931334451200}t^{10}+
\frac{4077796979}{643778150400}t^{11}+
\frac{8513133061}{9656672256000}t^{12}+ \frac{468626303}{4707627724800}t^{13}
+ \frac{26270857}{2897001676800}t^{14}+ \\
& \\
& \frac{124270847}{188305108992000}t^{15}+ \frac{2371609}{62768369664000}t^{16}
+ \frac{1182547}{711374856192000}t^{17}+ \frac{593}{10944228556800}t^{18}+
\frac{149789}{121645100408832000}t^{19}+\\
& \\
& \frac{2117}{121645100408832000}t^{20}
+ \frac{1}{8688935743488000}t^{21} \\
& \\ \hline

\end{array}
\]
}

We conclude with some remarks. First, it is natural to
ask whether one can derive our volume formula from a perturbation of
the Birkhoff polytope and then applying Lawrence's formula for simple
polytopes \cite{lawrence}. We found such a proof using a perturbation
suggested by B. Sturmfels, but the proof presented here yields more
results, for example, Corollaries \ref{twoells} and \ref{eqsigma} can
only be obtained this way.  Second, it is well known that Brion's and
Lawrence's formulas can be proved from the properties of
characteristic functions of polyhedra under polarity (see
e.g. Corollary 2.8 in \cite{BarviPom} or Theorem 3.2
\cite{beckhassesottile}). On the other hand P. Filliman
\cite{filliman, kuperberg} expressed the characteristic function any
convex polytope $P$ containing the origin as an alternating sum of
simplices that share supporting hyperplanes with $P$. The terms in the
alternating sum are given by a triangulation of the polar polytope of
$P$. Filliman's machinery yields in a limiting case Lawrence's volume
formulas. Different choices of triangulation of the polar of $P$ yield
different volume formulas for $P$. Using Filliman's duality
G. Kuperberg found (unpublished) other special volume formulas
that follow from pulling triangulations of the dual of $B_n$.

\section*{\bf Acknowledgments} 
We are truly grateful to Richard P. Stanley for his
encouragement and support. In fact it is because Prof. Stanley that
this project started, after he alerted us that there was a really nice
pattern in the data presented in Table 2 of \cite{latte}.  We are also
grateful to Bernd Sturmfels and G\"unter Ziegler for useful
conversations that eventually led us to discover the combinatorial
statement proved in Theorem \ref{main}. We have also received useful
suggestions and references from Alexander Barvinok, Matthias Beck, Greg Kuperberg,
Alexander Postnikov, and Peter Huggins.

\end{document}